\documentclass[11pt,reqno,a4pape]{amsart}
\usepackage{times,mathrsfs,color,comment}
\usepackage{amssymb,amsmath,bbm}
\usepackage[dvipsnames]{xcolor}

\numberwithin{equation}{section}

\newtheorem{thm}{Theorem}[section]
\newtheorem{cor}[thm]{Corollary}
\newtheorem{lem}[thm]{Lemma}

\newtheorem{prop}[thm]{Proposition}

\theoremstyle{remark}
\newtheorem{rem}[thm]{\bf Remark}

\newcommand{\ip}[2]{\langle #1,#2\rangle}
\newcommand{\triplenorm}[1]{\left|\!\left|\!\left|#1\right|\!\right|\!\right|}

\newcommand{\re}{{\rm Re\,}}
\newcommand{\C}{\mathbb{C}}

\newcommand{\T}{\mathbb T}

\newcommand{\F}{\mathcal F}
\newcommand{\norm}[1]{\left\lVert #1\right\rVert}
\newcommand{\abs}[1]{\left\lvert #1\right\rvert}
\newcommand{\Mat}{\mathscr M_N}
\newcommand{\w}{{\bf w}}
\newcommand{\cO}{\mathcal O}

\DeclareRobustCommand{\rchi}{{\mathpalette\irchi\relax}}
\newcommand{\irchi}[2]{\raisebox{\depth}{$#1\chi$}}

\title[The rigidity phenomenon for WCOs]{The rigidity phenomenon for weighted composition operators on the Fock spaces}

\date{November 19, 2026}

\author[H. Cho]{Hong Rae Cho}
\address{(Cho) Department of Mathematics, Quantum Science Technology Center, and Institute of Mathematical Sciences, Pusan National University, Pusan 46241, Republic of KOREA}
\email{chohr@pusan.ac.kr}

\author[B. Choe]{Boo Rim Choe}
\address{(Choe) Department of Mathematics, Korea University, Seoul 02841, Republic of KOREA}
\email{cbr@korea.ac.kr}

\author[H. Koo]{Hyungwoon Koo}
\address{(Koo) Department of Mathematics, Korea University, Seoul 02841, Republic of KOREA}
\email{koohw@korea.ac.kr}

\author[P. Tien]{Pham Trong Tien}
\address{(Tien) Faculty of Mathematics, Mechanics and Informatics, VNU University of Science, Vietnam National
University, Hanoi, Vietnam}
\email{phamtien@vnu.edu.vn}


\thanks{H. Cho was supported by a 2-Year Research Grant of Pusan National University. 
B. Choe was supported by NRF(RS-2026-25476546) of Korea.  H. Koo was supported by NRF(RS-2025-23324102) of Korea. Pham Trong Tien was funded by Vietnam National Foundation for Science and Technology Development (NAFOSTED) under grant number 101.02-2025.08.}

\begin{document}
\maketitle

\begin{abstract}
We establish the rigidity of canonical sums of any finite number of weighted composition operators acting from a Fock space $\mathcal{F}^p$ into another $\mathcal{F}^q$ for all possible combinations of  \(p,q\in (0, \infty]\). Specifically, we show that such a sum is bounded (compact, resp.) if and only if its summands are individually bounded (compact, resp.). In the case \(p\le q\), we present two independent approaches -- one analytic and one largely algebraic -- based on some common preliminary rigidity results. As a natural byproduct of our approach, we also obtain the Reproducing Kernel Thesis.
\end{abstract}

\section{Introduction}

Let $\cO(\C)$ be the space of all entire functions on the complex plane $\C$. Fix $\alpha>0$ throughout the paper.  

For $0<p<\infty$, denote by $\F^p$ the Fock space of all functions $f\in\cO(\C)$ such that
$$
\triplenorm{f}_p:=\left\{\frac{p\alpha}{2\pi}\int_{\C}\left|f(z)e^{-\frac{\alpha}{2}|z|^2}\right|^p \,dA(z)\right\}^{1/p}<\infty,
$$
where $dA$ denotes the area measure on $\C$. The triple-bar notation is adopted to distinguish the Fock (quasi-)norm from the standard Lebesgue norm. Meanwhile, the factor $\frac{p\alpha}{2\pi}$ is chosen solely for normalization purposes. When $p=\infty$, denote by $\F^\infty$ the space of all $f\in\cO(\C)$ such that
$$
\triplenorm{f}_\infty:=\sup_{z\in\C}|f(z)|e^{-\frac{\alpha}{2}|z|^2}<\infty.
$$
As is well known, $\F^p$ is a Banach space for $1\le p\le\infty$ and a complete metric space for $0<p<1$ with the metric 
$(f, g)\mapsto \triplenorm{f-g}_p^p$. The reference \cite{Zhu} is entirely devoted to these spaces and their operators. 

By an {\em inducing pair} $(\psi, \varphi)$ we mean an ordered pair $(\psi, \varphi)$ of {\em symbol} $\varphi\in \cO(\C)$ and {\em weight} $\psi\in \mathcal O(\C)$ with $\psi\not\equiv 0$.
Associated with an inducing pair $(\psi, \varphi)$ is the weighted composition operator (WCO) $W_{\psi, \varphi}$  defined by 
$$
W_{\psi, \varphi}f:= \psi \cdot (f\circ \varphi)
$$
for $f\in \mathcal O(\C)$. So, when $\psi\equiv1$, $W_{1, \varphi}$ is the standard composition operator (CO) acting on $\mathcal O(\C)$. 

Compared to the study of composition operators on Bergman or Hardy spaces (see classical references \cite{CM, Sh} for earlier studies), 
research on Fock spaces has a relatively short history. The study of COs on Fock spaces was initiated in 2003 by Carswell, MacCluer, and Schuster \cite{CMS}, who established fundamental characterizations for their boundedness and compactness on the several-variable Hilbert-Fock space. 
Motivated by \cite{CMS}, Ueki \cite{U} obtained characterizations (in the one-variable setting) for boundedness and compactness of WCOs acting from $\F^p$ to $\F^\infty$. 
More interestingly, Choe, Izuchi, and Koo \cite{CIK} found the $2$-rigidity phenomenon on the several-variable Hilbert-Fock spaces, meaning that a non-degenerate linear sum of two COs is bounded  (compact, resp.) if and only if its summands are individually bounded (compact, resp.). In subsequent one-variable works \cite{KTT, TK1, TK2, Y} by Khoi, Thom, Tien, and Yang, the aforementioned results were extended  to WCOs in general settings, allowing the integrability exponents of the domain and target spaces to differ without any restrictions.

Attempts have been made to extend the aforementioned $2$-rigidity phenomenon in two distinct directions. 
The first three authors of the current paper \cite{CCK1} announced the rigidity phenomenon for the linear sum of any finite number of  COs on the several-variable Hilbert-Fock space.
However, a critical gap was subsequently discovered in their argument (specifically \cite[line 1 of p.1244]{CCK1}), restricting the validity of the proof to the case where symbols are all polynomial maps.
Meanwhile, although Lieu et al. \cite{LTT} extended the same $2$-rigidity phenomenon to the general several-variable setting where the integrability exponents of the domain and target spaces are allowed to differ arbitrarily, their approach is not applicable to WCOs. 

The primary goal of this paper is to establish the $N$-rigidity phenomenon for WCOs for any positive integer $N$ and across all possible combinations of integrability exponents.
Unfortunately, our current techniques are restricted to the one-variable setting, leaving the several-variable case open for future investigation.

To state our result, we introduce some notation and terminology to be used throughout the paper. We fix an arbitrary positive integer $N$--\,the number of summands\,--\,and set 
$$
J:=\{1, \dots, N\}
$$
for convenience. We also fix a collection $\{(\psi_j, \varphi_j)\}_{j\in J}$ of $N$ inducing pairs and set
\begin{align}\label{eq:operator}
T:= \sum_{j=1}^N W_{\psi_j, \varphi_j}.
\end{align}
When this sum is considered, we may combine weights associated with the same symbol and discard any resulting zero term. So, we always assume that $\{(\psi_j, \varphi_j)\}_{j\in J}$ is a {\em canonical} collection, meaning that 
$$
\varphi_i\not\equiv\varphi_j\quad\text{whenever}\quad i \ne j.
$$
Finally, for $0<q\le\infty$, we define the {\em kernel-norm $q$-transform} of $T$ by 
$$
\widetilde T_q(w):=\triplenorm{Tk_w}_q\ \text{(possibly $\infty$)},\qquad w\in\C,
$$
where  $k_w$ denotes the normalized reproducing kernel at $w\in\C$ with the property that 
\begin{align}\label{eq:normalized}
\triplenorm{k_w}_p\equiv 1\qquad (0<p\le\infty);
\end{align}
see Section \ref{sec:reproducing-kernel}.
While $\widetilde T_q$ is lower semicontinuous in general, it is continuous when $\widetilde T_q\in L^\infty$; see Section \ref{sec:transform}. Here, $L^s=L^s(\C, dA)$, $0<s\le\infty$, denotes the standard Lebesgue space over $\C$.
Unless otherwise emphasized, the notation specified above will be used throughout the paper without any further reference.

As in the earlier works mentioned above, our results are also divided into two cases, namely, the case $0<p\le q\le\infty$ and the remaining case.
In each case, as a side benefit of our argument, we also obtain a characterization -- the Reproducing Kernel Thesis -- by means of kernel-norm transforms.

It seems worth mentioning that Assertion (b) below is independent of the parameter $p$ by \eqref{eq:normalized}.
In what follows $C_0$ denotes the class of all functions continuous on $\C$ and vanishing at $\infty$.

\begin{thm}\label{thm:main1}
Let $0<p\le  q\le \infty$ and $\{(\psi_j,\varphi_j)\}_{j\in J}$ be a canonical collection of inducing pairs. Then the following assertions are equivalent:
\begin{enumerate}
\item[(a)] $T:\F^p\to\F^q$ is bounded (compact, resp.);
\item[(b)] $\widetilde T_q \in L^\infty$ ($C_0$, resp.);
\item[(c)] $W_{\psi_j, \varphi_j}:\F^p\to\F^q$ is bounded (compact, resp.) for each $j\in J$. 
\end{enumerate}
\end{thm}

In the remaining case, our result is the next theorem.

\begin{thm}\label{thm:main2}
Let $0<q< p\le \infty$ and $\{(\psi_j,\varphi_j)\}_{j\in J}$ be a canonical collection of inducing pairs. 
Put $\lambda:=\frac{q}{p}$ which is understood to be $0$ for $p=\infty$.
Then the following assertions are equivalent:
\begin{enumerate}
\item[(a)] $T:\F^p\to\F^q$ is bounded;
\item[(b)] $T:\F^p\to\F^q$ is compact;
\item[(c)] $\widetilde T_q \in L^{\frac{q}{1-\lambda}}$;
\item[(d)] $W_{\psi_j, \varphi_j}:\F^p\to\F^q$ is bounded for each $j\in J$;
\item[(e)] $W_{\psi_j, \varphi_j}:\F^p\to\F^q$ is compact for each $j\in J$. 
\end{enumerate}
\end{thm}

For the proof of Theorem \ref{thm:main1}, we present two independent approaches after establishing common preliminary rigidity results for both symbols and weights. The first is an analytic proof based on a slope-class isolation argument by means of appropriate differential operators, while the second is an algebraic approach grounded in uniform operator norm estimates for Gram matrices associated with finite systems of normalized reproducing kernels at separated points. Moreover, this algebraic approach plays the key role for our proof of Theorem \ref{thm:main2}. We believe each proof offers distinct merits and valuable perspectives. 
\medskip

Broadly speaking, the proofs of our two main theorems share a preliminary rigidity step and then proceed along different routes, depending on the relation between $p$ and $q$.
\begin{itemize}
\item \textbf{Preliminary rigidity (Section~3).} A differential annihilation and logarithmic averaging argument for exponential sums forces all symbols to be affine. When all symbols are constant, a Vandermonde-matrix argument forces all weights to belong to the target space.
\item \textbf{The range $0<p\le q\le\infty$ (Section~4).} We establish the boundedness and compactness of individual summands by controlling their gauge functions. Two independent approaches
are presented. The first isolates each slope-class operator by means of suitable differential operators (Lemma~4.2). The second uses uniform Gram-matrix estimates for normalized reproducing
kernels at separated points (Theorem~4.3).
\item \textbf{The range $0<q<p\le\infty$ (Section~5).} Building on the Gram-matrix approach, we combine lattice discretization, uniform dual-system estimates, a Khinchine-type randomization argument, and sequence-space duality to establish gauge integrability (Proposition~5.7). Together with Gaussian convolution estimates, this yields the kernel-norm characterization in Theorem~1.2.
\end{itemize}
\medskip
 
The paper is organized as follows:
\smallskip

In Section \ref{sec:prelim} we recall and collect some basic notions and related facts. 

In Section \ref{sec:pre-rigidity}, as a preliminary step toward our main results, we establish the rigidity properties of both symbols and weights. 

In Section \ref{sec:main1} we establish Theorem \ref{thm:main1}. Two independent proofs, one analytic and one largely  algebraic, are presented.

In Section \ref{sec:main2} we establish Theorem \ref{thm:main2}. The proof is achieved through a nontrivial application of the algebraic approach in the previous section.

\medskip
{\it Constants.} Throughout the paper we use the same letter $C$ to denote positive constants
which may vary at each occurrence but do not depend on the essential parameters.
The dependence of constant $C$ on specific parameters will be indicated when relevant.
For nonnegative quantities $X$ and $Y$ the notation $X\lesssim Y$  or $Y\gtrsim X$
means $X\le CY$ for some inessential constant $C$. Similarly, we write $X\approx Y$ if
both $X\lesssim Y$ and $Y\lesssim X$ hold.

\section{Preliminaries}\label{sec:prelim}

In this section, we collect basic results on key tools required for later sections.

\subsection{Compactness criterion}
Let $0<p, q\le \infty$ and $L$ be a sum of finitely many WCOs. Recall that 
$L:\F^p\to \F^q$ is said to be {\em compact} if it is bounded, meaning that 
$$
\|L\|_{\F^p\to \F^q}:=\sup_{\triplenorm{f}_p= 1}\triplenorm{Lf}_q<\infty,
$$
and if its image of each bounded set in $\F^p$ is relatively compact in $\F^q$.

We have a convenient compactness criterion asserting that $L$ is compact if and only if 
$\triplenorm{Lf_n}_q\to 0$ for any bounded sequence $\{f_n\}$ in $\F^p$ that converges to $0$ locally uniformly; see \cite[Lemma 1.2]{KTT} or \cite[Lemma 2.3]{TK2}.

\subsection{Reproducing kernels} \label{sec:reproducing-kernel}
When $p=2$, the space $\F^2$ is a Hilbert space with inner product 
$$
\ip{f}{g}:=\frac{\alpha}{\pi} \int_\C f(z)\overline{g(z)} e^{-\alpha|z|^2}\, dA(z).
$$
In fact  $\F^2$ is a reproducing kernel Hilbert space whose reproducing kernel at $w\in \C$ is given by  
$$
K_w(z):=e^{\alpha\overline{w} z}.
$$
It is easy to verify that $\triplenorm{K_w}_p= e^{\frac{\alpha}{2}|w|^2}$ for any $0<p\le \infty$. Hence, the $\F^p$-normalized reproducing kernel
\begin{align*}
k_w(z):=\frac{K_w(z)}{\triplenorm{K_w}_p}=e^{\frac{\alpha}{2}(2\overline  w z -|w|^2)}
\end{align*}
is independent of the parameter $p$ and its modulus admits a neat formula 
\begin{align}\label{eq:kernelsize}
|k_w(z)|=e^{\frac{\alpha}{2}(|z|^2-|z-w|^2)},
\end{align}
which yields \eqref{eq:normalized}. It is  clear from the above that $k_w\to 0$ locally uniformly as $|w|\to \infty$.

\subsection{Submean estimates}\label{sec:submean-estimate}
Let $f\in\cO(\C)$ and $w\in\C$. One may easily verify by the reproducing property and the Cauchy-Schwarz Inequality 
$$
|f(w)| e^{-\frac{\alpha}{2}|w|^2}\le \triplenorm{f}_2.
$$
This inequality is known to extend to all $0<q\le \infty$, i.e.,  
\begin{equation}\label{eq:point-eval}
|f(w)|e^{-\frac{\alpha}{2}|w|^2}\le \triplenorm{f}_q;
\end{equation}
see \cite[Corollary 2.8]{Zhu}. In addition, when $0<q<\infty$,
\begin{align}\label{eq:vanish-infty}
\lim_{|w|\to\infty}|f(w)|e^{-\frac{\alpha}{2}|w|^2}=0\quad\text{for}\quad f\in \F^q;
\end{align}
see \cite[Proposition 2.9]{Zhu}.

A local version of \eqref{eq:point-eval}
is also available. Namely, denoting by
$$
D_r(w):=\{z\in\C: |z-w|<r\}
$$
the disk with center $w\in\C$ and radius $r>0$, we have
\begin{align}\label{eq:local}
|f(w)| e^{-\frac{\alpha}{2}|w|^2}\le C  \left\{\int_{D_r(w)} \left|f(z)e^{-\frac{\alpha}{2}|z|^2}\right|^q\, dA(z) \right\}^{1/q}
\end{align}
for some constant $C=C(\alpha, q, r)>0$; see \cite[Lemma 3.1]{CCK2} or \cite[Lemma 2.32]{Zhu}.

Applying \eqref{eq:local} to the entire function  $w\mapsto \overline {TK_w(z)}$ (with $z$ arbitrary but fixed) gives
\begin{align*}
|Tk_w(z)|^q\lesssim \int_{D_r(w)} \left| Tk_\zeta(z) \right|^q \, dA(\zeta)
\end{align*}
Integrating against the relevant Gaussian measure and interchanging the order of integrations yields
\begin{align}\label{eq:submean-transform}
\widetilde T^q_q(w) \le C \int_{D_r(w)} \widetilde T^q_q(\zeta)\,dA(\zeta)
\end{align}
for some constant $C=C(\alpha, q, r)>0$.

\subsection{Single operators}\label{sec:single-operator}
Associated with an inducing pair $(\psi,\varphi)$ is the {\em gauge} function 
\begin{equation*}
M_{(\psi,\varphi)}(z):=|\psi(z)|e^{\frac\alpha2(|\varphi(z)|^2-|z|^2)}.
\end{equation*}
This gauge function naturally comes into play in the study of boundedness and compactness of $W_{\psi,\varphi}$ acting from a Fock space into another. 
First of all, when $p < \infty$ or $q < \infty$, the boundedness of the operator $W_{\psi,\varphi}:\mathcal F^p\to\mathcal F^q$ implies 
\begin{align*}
M_{(\psi,\varphi)}\in L^\infty,
\end{align*}
which in turn implies that $\varphi$ is affine and $|\varphi'(0)|\le 1$, i.e., 
\begin{align}\label{eq:affine}
\varphi(z)=az+b\ \ \text{with} \ \ |a|\le1;
\end{align}
see \cite[Proposition~3.1]{TK2} and \cite[Propositions~2.1 and 3.1]{KTT}.

Conversely, for an affine symbol $\varphi$ as in \eqref{eq:affine}, the boundedness
and compactness of $W_{\psi,\varphi}$ are characterized as in the next proposition. 

\begin{prop} \label{prop:single-wco}
Let $0<p,q\le\infty$ and $(\psi, \varphi)$ be an inducing pair with $\varphi$ as in \eqref{eq:affine}.
Then the following assertions hold:
\begin{enumerate}
\item[(a)] Let $a=0$. Then $W_{\psi,\varphi}:\F^p\to \F^q$ is bounded if and only if
$\psi\in\mathcal F^q$. In this case, it has rank one, hence is compact; 
\item[(b)] Let $a\ne 0$ and $0<p\le q\le\infty$. Then $W_{\psi,\varphi}:\mathcal F^p\to\mathcal F^q$ is bounded (compact, resp.)  if and only if $M_{(\psi,\varphi)}\in L^\infty$ ($C_0$, resp.);
\item[(c)] Let $a\ne 0$ and $0<q<p\le\infty$. Put $\lambda:=\frac{q}{p}$ which is understood to be $0$ for $p=\infty$. Then $W_{\psi,\varphi}:\mathcal F^p\to\mathcal F^q$ is bounded/compact if and only if $M_{(\psi,\varphi)} \in L^{\frac{q}{1-\lambda}}$.
\end{enumerate}
\end{prop}

\begin{proof}
Note that point evaluations are bounded on $\F^p$ by \eqref{eq:point-eval}. Also, note $\psi=W_{\psi,\varphi}1$. Thus (a) holds. 

To prove (b), we first consider the case $p=q=\infty$. In this case, we note 
\begin{align*}
|W_{\psi,\varphi}f(z)|e^{-\frac{\alpha}{2}|z|^2} \le M_{(\psi,\varphi)}(z)
\le \triplenorm{W_{\psi,\varphi}k_{\varphi(z)}}_\infty
\end{align*}
for all $z\in\C$ and $f\in \F^\infty$ with $\triplenorm{f}_\infty\le 1$. Also, note $\triplenorm{k_{\varphi(z)}}_\infty=1$ for all $z\in\C$ and $k_{\varphi(z)}\to 0$ locally uniformly as $|z|\to\infty$. Thus one may easily verify that (b) holds in the case $p=q=\infty$.

Now, consider the case $p<\infty$ or $q<\infty$. In this case, under the additional assumption $\psi\in \F^q$, (b) and (c) are proved in \cite[Theorems 3.3 and 3.4]{TK2} and \cite[Theorems 2.3 and 3.3]{KTT}. So, in order to prove (b), we assume $M_{(\psi,\varphi)}\in L^\infty$ and show $\psi\in \F^q$. For $q=\infty$, we note 
\begin{align}\label{eq:Midentity}
|\psi(z)| e^{-\frac{\alpha }{2}|z|^2}=M_{(\psi,\varphi)}(z)e^{-\frac{\alpha}{2}|az+b|^2}
\end{align}
and thus 
\begin{align*}
\triplenorm{\psi}_\infty\le \norm{M_{(\psi,\varphi)}}_{L^\infty}.
\end{align*}
Also, for $q<\infty$, we note
\begin{align*}
\triplenorm{\psi}^q_q  \le \frac{q\alpha}{2\pi}\norm{M_{(\psi,\varphi)}}^q_{L^\infty} \int_{\C}e^{-\frac{q\alpha}{2}|az+b|^2}\,dA(z) <\infty
\end{align*}
where the last inequality holds by an elementary change-of-variable. Thus (b) holds. 

Now, in order to prove (c), we assume $M_{(\psi,\varphi)} \in L^{s}$ where $s:=\frac{q}{1-\lambda}$ and show $\psi\in\F^q$. 
When $p<\infty$, using \eqref{eq:Midentity} and applying H\"older's Inequality based on the exponent relation $\frac{q}{s}+ \lambda=1$, we 
\begin{align*}
\begin{aligned}
\triplenorm{\psi}^q_q  
&\le \frac{q\alpha}{2\pi}\norm{M_{(\psi,\varphi)}}^q_{L^{s}} \left\{\int_{\C}e^{-\frac{p\alpha}{2}|az+b|^2}\,dA(z)\right\}^{\lambda}<\infty,
\end{aligned}
\end{align*}
as desired. This remains valid when $p=\infty$ with $s=q$ and $\lambda=0$. The proof is complete.
\end{proof}

Once Theorems \ref{thm:main1} and \ref{thm:main2} are established, the gauge conditions in Proposition \ref{prop:single-wco} can be replaced by their corresponding counterparts for kernel-norm transforms. Of course, such results can also be derived directly for single operators; in fact, the case $0<q<p<\infty$ is implicit in the proof of \cite[Theorem 3.3]{TK2}, and  other cases are handled in a standard way (see Remark \ref{remark} at the end of the paper).

\subsection{Gram matrix: Abstract setting}\label{sec:Gram=abstract}
Denote by $\Mat$ the space of all $N\times N$ complex matrices. We mean by $L=(a_{ij})\in\Mat$ that the $ij$-entry of $L$ is $a_{ij}$. 
We will identify $L\in\Mat$ with the linear operator $L:\C^N\to \C^N$. We will use the notation $\norm{L}$  
for the operator norm of $L$ regarded as an operator acting on $\C^N$.

Let $\mathcal H$ be an $N$-dimensional Hilbert space with inner product $\ip{\cdot}{\cdot}_{\mathcal H}$. Let ${\bf x}=(x_j)_{j=1}^N$ with $x_j\in\mathcal H$ for each $j$. 
Associated with ${\bf x}$ is the Gram matrix 
$$
G_{\bf x}=\bigl(\ip{x_j}{x_i}_{\mathcal H}\bigr)\in \Mat.
$$
Clearly, $G_{\bf x}$ is self-adjoint. It is elementary to verify that $G_{\bf x}$ is invertible if and only if ${\bf x}$ is an $N$-tuple of linearly independent vectors.
Furthermore, $G_{\bf x}$ is well known to be positive semi-definite. In fact, denoting by $E_{\bf x}:\C^N\to \mathcal H$ the linear operator defined by $E_{\bf x}{\bf e}_j=x_j$ for each $j\in J$ where $({\bf e}_j)_{j\in J}$ denotes the standard basis for $\C^N$, one may easily verify 
\begin{align}\label{eq:positive}
G_{\bf x}= E_{\bf x}^* E_{\bf x}.
\end{align}
In particular, 
\begin{align}\label{eq:eigen-max}
\norm{G_{\bf x}}= \lambda_{\rm max}(G_{\bf x})
\end{align}
where $\lambda_{\rm max}$ denotes the largest eigenvalue; see \cite[Example~5.6.6]{HJ}.
 
Now assume that ${\bf x}$ is an $N$-tuple of linearly independent vectors, and hence that $\bf x$ is a basis for $\mathcal H$. Let $\beta=(\beta_j)_{j\in J}$ be the dual basis of $\bf x$, meaning that $\beta$ is the basis of $\mathcal H$ uniquely determined by the biorthogonality 
\begin{equation}\label{eq:biorth}
\ip{\beta_j}{x_i}_{\mathcal H}=\delta_{ij}\qquad(i, j\in J)
\end{equation}
where $\delta_{ij}$ denotes the Kronecker delta. To describe $\beta$ more explicitly, let 
$$
G^{-1}_{\bf x}=:\bigl( y_{ij}\bigr).
$$
Note that the equality $G_{\bf x}G_{\bf x}^{-1}=I_N$ is precisely the same as the biorthogonality
$$
\delta_{ij}= \sum_{\ell=1}^N \ip{x_\ell }{x_i}_{\mathcal H} y_{\ell j} =\left \langle{\sum_{\ell=1}^N  y_{\ell j} x_\ell }, {x_i}\right\rangle_{\mathcal H}\qquad(i,j\in J).
$$
Consequently, the representation of $\beta_j$ with respect to the basis $\bf x$ is
\begin{align*}
\beta_j= \sum_{\ell=1}^N  y_{\ell j} x_\ell
\end{align*}
for each $j$. Combined with \eqref{eq:biorth}, this in turn yields $\ip{\beta_j}{\beta_i}_{\mathcal H}= y_{ij}$ and thus 
\begin{align}\label{eq:inverse}
G_\beta= G_{\bf x}^{-1}.
\end{align}

\subsection{Gram matrix: Current setting}\label{sec:Gram}
We now set some notation and observe some basic properties for the setting of the current paper. 
We will use the notation 
$$
T{\bf f}:=(Tf_j)_{j\in J}
$$
for an $N$-tuple ${\bf f}=(f_j)_{j\in J}$ of entire functions.

Consider pairwise distinct complex numbers $w_1, \dots, w_N$. Note that $k_{w_1}, \dots, k_{w_N}$ are linearly independent.  
Put   
$$
\w:=(w_j)_{j\in J} \in\C^N\quad \text{and}\quad k_{\bf w}:= (k_{w_j})_{j\in J}.
$$
Since $k_{\w}$ is an $N$-tuple of linearly independent functions, it is a basis for the $N$-dimensional linear subspace $\operatorname{span} k_{\w}$ of $\F^2$. 
Let 
$$
\beta^\w=(\beta^\w_j)_{j\in J}
$$ 
be the dual basis of $k_\w$. We note from \eqref{eq:inverse}
\begin{align}\label{eq:Gram-inverse}
G_{\beta^\w}=G_{k_\w}^{-1}.
\end{align}
Now, fix $z\in\C$ and consider the case where
$$
\w=(\varphi_j(z))_{j\in J}
$$
with $\varphi_j(z)$'s pairwise distinct. 
By the reproducing property and the duality we have
$$
\beta^\w_{j} \bigl(\varphi_i(z)\bigr)= \ip{\beta^\w_{j}} {K_{\varphi_i(z)}} =\delta_{ij}e^{\frac{\alpha}{2} |\varphi_i(z)|^2}\qquad(i, j\in J)
$$
and thus 
\begin{align}\label{eq:duality}
T \beta^{\w} (z)
&= \left(\psi_j(z) e^{\frac{\alpha}{2}|\varphi_j(z)|^2}\right)_{j\in J}.
\end{align}
Accordingly, introducing the combined gauge
\begin{align}\label{eq:gaugenorm}
\mathcal M(z):=\left\{\sum_{j=1}^N M^2_{(\psi_j,\varphi_j)}(z) \right\}^{1/2},
\end{align}
we note 
\begin{align}\label{eq:M-equality}
\mathcal M(z)=|T\beta^\w(z)|e^{-\frac{\alpha}{2}|z|^2}.
\end{align}
Here, the same notation $|\cdot|$, which is being used for the modulus of a complex number, is used for the Euclidean norm on $\C^N$.  
The meaning of these notations should be clear from the context. 
Furthermore, we see from \eqref{eq:duality} that
\begin{align*}
T \beta^\w_{j} (z)
&=\sum_{i=1}^N \psi_i(z) \delta_{ij}e^{\frac{\alpha}{2}|w_i|^2}\\
&=\sum_{i=1}^N \psi_i(z) \left\{\sum_{\ell=1}^N \overline{\ip{\beta^\w_{\ell}}{\beta^\w_{j}}\ip{k_{w_i}} {k_{w_\ell}}}\right\}e^{\frac{\alpha}{2}|w_i|^2}\qquad\text{by \eqref{eq:Gram-inverse}}\\
&=\sum_{\ell=1}^N \overline{\ip{\beta^\w_\ell}{\beta^\w_j}}\sum_{i=1}^N \psi_i(z) \ip{k_{w_\ell}}{K_{w_i}}\\
&=\sum_{\ell=1}^N \overline{\ip{\beta^\w_\ell}{\beta^\w_j}}  Tk_{w_\ell}(z),
\end{align*}
which means that the $j$-th entry of $\overline{G_{\beta^\w}}[Tk_\w (z)]$ is $T\beta^\w_j(z)$ for each $j$. Accordingly, 
\begin{align}\label{eq:duality2}
\overline{G_{\beta^\w}}[Tk_\w (z)]= T\beta^\w(z).
\end{align}
The introduction of the Gram matrix and the dual basis is motivated by equalities \eqref{eq:Gram-inverse}, \eqref{eq:duality} and \eqref{eq:duality2}. 

\subsection{Kernel-norm transform: Continuity}\label{sec:transform}
Let $0<q\le \infty$. We first notice that $\widetilde T_q$ is lower semicontinuous. 
One may verify it by Fatou's Lemma when $q<\infty$, and by the fact that $\widetilde T_\infty$ is the supremum of positive continuous functions.
In particular, the supremum and the essential supremum are the same for $\widetilde T_q$, i.e.,  
\begin{align*}
\sup_{w\in\C} \widetilde T_q(w) = \|\widetilde T_q\|_{L^\infty}.
\end{align*}
This will be tacitly used in what follows.
We also notice that the following implication holds: 
\begin{align}\label{eq:transform-continuity}
\text{$\widetilde T_q\in L^\infty$}\implies \text{$\widetilde T_q:\C\to [0,\infty)$ is continuous}.
\end{align}
({\em Sketch of Proof}): 
Assume $\widetilde T_q\in L^\infty$. Consider arbitrary $z, w_0\in\C$ and let $|w-w_0|\le 1$. Applying the Cauchy Integral Formula to the entire function $\zeta\mapsto \overline {TK_{\zeta}(z)}$ along the circle $|\zeta-w_0|=2$, we obtain
\begin{align}
|T K_{w}(z)-T K_{w_0}(z)| 
&\le  \frac{1}{2\pi}\int_{|\zeta-w_0|=2} \left|\frac{TK_{\zeta}(z)(w-w_0)}{(\zeta-w)(\zeta-w_0)}\right|\, |d\zeta|\nonumber\\
&\le 2 |w-w_0|\left(\sup_{|\zeta-w_0|=2} |T K_{\zeta}(z)|\right)\label{eq:pre-cont}.
\end{align}
Moreover, when $0<q<\infty$, subharmonicity in $\zeta$-variable over the closed disk $\overline D_1(\zeta)$ yields 
\begin{align*}
\sup_{|\zeta-w_0|=2}|TK_\zeta(z)|^q  \le \frac{1}{\pi}\int_{D_3(w_0)} |TK_\eta(z)|^q\, dA(\eta).
\end{align*}
Combined with \eqref{eq:pre-cont}, this yields 
\begin{align*}
|T K_{w}(z)-T K_{w_0}(z)|^q\lesssim |w-w_0|^q\int_{D_3(w_0)} |TK_\eta(z)|^q\, dA(\eta)
\end{align*}
with hidden constant depending only on $q$. Now, integrating against the relevant Gaussian measure and interchanging the order of integrations, we obtain 
\begin{align*}
\triplenorm{TK_{w} -T K_{w_0}}_q \lesssim |w- w_0|\bigl\|\widetilde T_q\bigr\|_{L^\infty}.
\end{align*}
Note that this remains valid for $q=\infty$ directly by \eqref{eq:pre-cont}. 
This implies that $w\mapsto TK_w$ is a continuous function from $\C$ into $\F^q$. 
As a consequence, since 
$$
\widetilde T_q (w)= e^{-\frac{\alpha}{2}|w|^2} \triplenorm{TK_w}_q,
$$
we see that $\widetilde T_q:\C\to [0, \infty)$ is continuous. $\Box$.

The implication  
\begin{align}\label{eq:vanish}
\text{$T:\F^p\to \F^q$ is compact  $\implies \widetilde T_q\in C_0$}
\end{align}
is now clear by \eqref{eq:transform-continuity} and the compactness criterion.

\section{The Symbol/Weight Rigidity}\label{sec:pre-rigidity}

In this section, as a preliminary step toward separation of summands, we first show that the operator in \eqref{eq:operator} can be bounded only when all the symbols are affine.
Once all symbols are established to be affine, note from Proposition \ref{prop:single-wco} that constant-symbol summands require separate treatment.
Taking this into account, we show that their sum can be bounded only when all the weights belong to the target space.

We begin with the following exponential-sum rigidity lemma, which is a bit more general than needed.

\begin{lem}\label{lem:rigidity}
Let $\{(\psi_j,\varphi_j)\}_{j\in J}$ be a canonical collection of inducing pairs. Suppose that there are constants $c_1>0$, $c_2\in(0,\alpha)$, and $C>0$ such that
\begin{equation}\label{eq:rigidity-hyp}
\left|\sum_{j=1}^N \psi_j(z)e^{\alpha w\varphi_j(z)}\right|\le C e^{c_1|z|^2+c_2|w|^2}
\end{equation}
for all $z,w\in\C$. Then each $\varphi_j$ is affine.
\end{lem}

\begin{proof}
Fix $z\in\C$ and put
\begin{align*}
F_z(w):=\sum_{j=1}^N\psi_j(z)e^{\alpha w\varphi_j(z)}\in\cO(\C).
\end{align*}
Given $k\in J$, set 
\begin{align*}
L_k:=\prod_{\ell\ne k}\bigl(\partial_w-\alpha\varphi_\ell(z)\bigr)\quad\text{and}\quad 
P_k(z):=\prod_{\ell\ne k}\bigl(\varphi_k(z)-\varphi_\ell(z)\bigr)
\end{align*}
where $\partial_w$ denotes the complex differential operator in the $w$-variable. Note 
\begin{equation*}
L_kF_z(w)
=\alpha^{N-1}\psi_k(z)P_k(z)e^{\alpha w\varphi_k(z)}.
\end{equation*}
On the other hand, we note
$$
L_k F_z(w)= F_z^{(N-1)}(w) + \sum_{m=1}^{N-1} (-1)^m \alpha^m e_m\bigl(\check\Phi_k(z)\bigr) F_z^{(N-m-1)}(w)
$$
where $\check\Phi_k(z)=\bigl(\varphi_m(z)\bigr)_{m\ne k}$ and  $e_m$ is the $m$-th elementary symmetric polynomial in $(N-1)$ variables. 
Note 
$$
\left| e_m\bigl(\check\Phi_k(z)\bigr)\right|\lesssim M^m(z)\quad\text{where}\quad M(z):=\max_{j\in J}|\varphi_j(z)|.
$$
Accordingly, we have 
\begin{align}
|\psi_k(z)P_k(z)e^{\alpha w\varphi_k(z)}|
&=\alpha^{1-N}|L_k F_z(w)|\nonumber\\
&\lesssim |F_z^{(N-1)}(w)|+  \sum_{m=1}^{N-1} |M^m(z)| |F_z^{(N-m-1)}(w)|. \label{eq:F-derivative-growth}
\end{align}
Meanwhile, applying the Cauchy Estimates over the disk $D_{(1+|w|)^{-1}}(w)$, we obtain 
\begin{align*}
\frac{|F_z^{(j)}(w)|}{j!} 
&\le (1+|w|)^j \left(\sup_{|\zeta-w|=\frac{1}{1+|w|}} |F_z(\zeta)|\right)\\
&\lesssim (1+|w|)^j  e^{c_1|z|^2 + c_2|w|^2}\quad \text{by} \  \eqref{eq:rigidity-hyp}
\end{align*} 
for each $j$; note $|\zeta|^2 \le |w|^2 + 3$  for the last estimate. Combining this with \eqref{eq:F-derivative-growth}, we obtain
\begin{align}\label{eq:LkFest}
|\psi_k(z)P_k(z)e^{\alpha w\varphi_k(z)}|
\lesssim
(1+M(z)+|w|)^{N-1}e^{c_1|z|^2+c_2|w|^2}
\end{align}
with implicit constant independent of $z$ and $w$. Now, choosing $w=\overline{\varphi_k(z)}$, we obtain 
\begin{align*}
|\psi_k(z)P_k(z)|
\lesssim
\bigl(1+M(z)\bigr)^{N-1}
\exp\!\left\{c_1|z|^2-(\alpha-c_2)|\varphi_k(z)|^2\right\}.
\end{align*}
Recall that this estimate holds for each $k\in J$. So, taking the product over all $k$ and using 
$\sum_{k=1}^N |\varphi_k(z)|^2 \ge M^2(z)$, we obtain  
\begin{equation*}
\prod_{k=1}^N|\psi_k(z)P_k(z)|
\lesssim
(1+M(z))^{N(N-1)}
\exp\!\left\{Nc_1|z|^2-(\alpha-c_2)M^2(z)\right\}.
\end{equation*}
Since a fixed polynomial is dominated by any positive Gaussian, this yields 
\begin{equation}\label{eq:rigidity-product}
e^{\frac{\alpha-c_2}{2} M^2(z)}
\prod_{k=1}^N|\psi_k(z)P_k(z)|
\lesssim e^{Nc_1|z|^2}
\end{equation}
with hidden constant independent of $z$. 
Although the affine property of the symbols can now be derived via a related logarithmic averaging argument, as in \cite[Proposition 2.1]{L}, we include the details below to accommodate the more general growth condition \eqref{eq:rigidity-hyp}.

Since $\{(\psi_j,\varphi_j)\}_{j\in J}$ is a canonical collection of inducing pairs by assumption, we see that 
\begin{align*}
H(z):=\prod_{k=1}^N\psi_k(z)P_k(z)
\end{align*}
is a nonzero entire function. Write $H(z)=z^nH_0(z)$ where $n$ is a nonnegative integer and $H_0\in\cO(\C)$ with $H_0(0)\ne0$. From \eqref{eq:rigidity-product} we can find a constant $C$ independent of $z$ such that
\begin{align}\label{eq:pre-average}
\frac{\alpha-c_2}{2}|\varphi_j(z)|^2 +n\log|z|+\log|H_0(z)|-Nc_1|z|^2\le C
\end{align}
for each $j\in J$. 

Fix $j$ and $r>0$. Averaging over the circle $|z|=r$, we have by subharmonicity 
\begin{align*}
\int_\T\log|H_0(r\zeta)|\,d\sigma(\zeta)\ge\log|H_0(0)|
\end{align*}
where $\sigma$ is the normalized arc-length measure on the unit circle $\T$, and by Parseval's Identity
\begin{align*}
\int_{\T}|\varphi_j(r\zeta)|^2\,d\sigma(\zeta)=\sum_{\ell=0}^\infty
\left|\frac{\varphi_j^{(\ell)}(0)}{\ell!}\right|^2r^{2\ell}.
\end{align*}
Consequently, averaging both sides of \eqref{eq:pre-average} over the circle $|z|=r$, we obtain 
\begin{equation*}
\frac{\alpha-c_2}{2}
\sum_{\ell=0}^\infty
\left|\frac{\varphi_j^{(\ell)}(0)}{\ell!}\right|^2r^{2\ell}+n\log r+\log|H_0(0)|-Nc_1 r^2\le C.
\end{equation*}
Since this is valid for all $r>0$, dividing both sides of the above by $r^2$ and then taking the limit $r\to\infty$, we conclude 
$\varphi_j^{(\ell)}(0)=0$ for all $\ell\ge2$ and hence that $\varphi_j$ is affine.
The proof is complete. 
\end{proof}

Now, applying Lemma \ref{lem:rigidity}, we show that boundedness of $T$ along the normalized reproducing kernels yields the symbol rigidity as in the next proposition.  

\begin{prop}\label{prop:affine}
Let $0<q\le\infty$ and  $\{(\psi_j,\varphi_j)\}_{j\in J}$ be a canonical collection of inducing pairs. If $\widetilde T_q\in L^\infty$, 
then each $\varphi_j$ is affine.
\end{prop}

\begin{proof}
Let $z, w\in\C$. Applying $T$ to $k_{\overline w}$, we obtain by \eqref{eq:point-eval} 
\begin{align*}
|Tk_{\overline w}(z)|e^{-\frac{\alpha}{2}|z|^2}
\le \widetilde T_q(w)\le \bigl\|\widetilde T_q\bigr\|_{L^\infty}.
\end{align*}
On the other hand, we have 
\begin{align*}
Tk_{\overline w}(z)=e^{-\frac{\alpha}{2}|w|^2}\sum_{j=1}^N\psi_j(z)e^{\alpha w\varphi_j(z)}.
\end{align*}
Accordingly, we get
\begin{align}\label{eq:Tbdd}
\left|\sum_{j=1}^N\psi_j(z)e^{\alpha w\varphi_j(z)}\right|
\le e^{\frac\alpha2(|z|^2+|w|^2)} \bigl\|\widetilde T_q\bigr\|_{L^\infty}.
\end{align}
Lemma \ref{lem:rigidity} now applies with $c_1=c_2=\frac{\alpha}{2}$. The proof is complete.
\end{proof}

Having established symbol rigidity, we now turn to weight rigidity for sums of constant-symbol summands, as stated in the next proposition.
Note from Proposition \ref{prop:single-wco}(a) that the conclusion below can be rephrased as ``each summand is compact from $\F^p$ to $\F^q$ for any $0<p\le\infty$".

\begin{prop}\label{prop:weight-rigid}
Let $0<q\le\infty$ and  $\{(\psi_j,\varphi_j)\}_{j\in J}$ be a canonical collection of inducing pairs. Assume that all symbols are constant. 
If $\widetilde T_q\in L^\infty$, then 
$$
\psi_j\in \F^q
$$ 
for each $j\in J$. 
\end{prop}

\begin{proof}
Put $b_j:\equiv \varphi_j$ for $j\in J$. 
By canonicality, $b_1,\ldots,b_N$ are pairwise distinct. Note 
$$
Tk_w =e^{-\frac{\alpha}{2}|w|^2}\sum_{j=1}^{N}\psi_j e^{\alpha \overline w b_j} \in \mathcal F^q
$$
for all $w\in\C$. In particular, choosing $\epsilon>0$ so small that
$$
e^{\alpha\epsilon b_1}, \ldots, e^{\alpha\epsilon b_N}
$$
are pairwise distinct, we have 
\begin{align}\label{eq:Tkme}
e^{\frac{\alpha}{2}m^2\epsilon^2}Tk_{m\epsilon}=\sum_{j=1}^{N}\psi_j e^{m\alpha\epsilon b_j}\in \mathcal F^q
\end{align}
for $m=0,\ldots,N-1$.

Now, using the invertible Vandermonde matrix 
\begin{align}\label{eq:V}
V:=\bigl( e^{(i-1)\alpha\epsilon b_j} \bigr),
\end{align}
we can rephrase \eqref{eq:Tkme} as
\begin{align*}
\begin{pmatrix}  \psi_1 \\
\psi_2\\
 \vdots  \\
 \psi_N \\
 \end{pmatrix}
= V^{-1}
 \begin{pmatrix}  T k_0 \\
 e^{\frac{\alpha}{2}\epsilon^2} Tk_{\epsilon}\\
 \vdots  \\
e^{\frac{\alpha}{2}(N-1)^2\epsilon^2} Tk_{(N-1)\epsilon} \\
 \end{pmatrix}.
\end{align*}
With this relation established, one may now easily conclude the proposition. The proof is complete.
\end{proof}

Before closing the section, we introduce some notation suggested by Propositions \ref{prop:single-wco} and \ref{prop:affine}.

\subsection*{\quad\bf Notation (N)} Let $0<q\le \infty$ and assume $\widetilde T_q\in L^\infty$. By Proposition \ref{prop:affine}, all the symbols are affine. Pick the slope $a_j$ and the constant term $b_j$ of $\varphi_j$ so that    
\begin{align*}
\varphi_j(z)= a_j z + b_j\qquad (j\in J).
\end{align*}
By Proposition \ref{prop:single-wco}, constant-symbol summands will always be treated separately in the proofs of Theorems \ref{thm:main1} and \ref{thm:main2}. To this end, we fix some notation for use in later sections. Put
$$
J_0:=\{j\in J:a_j=0\}\quad\text{and}\quad J_1:=J\setminus J_0.
$$
Also, put  
\begin{align*}
T_0:=\sum_{i\in J_0} W_{\psi_i,\varphi_i}\quad\text{and}\quad T_1:= \sum_{i\in J_1} W_{\psi_i,\varphi_i}
\end{align*}
so that $T=T_0+T_1$. Finally, for $m=0, 1$, we set 
\begin{align*}
\mathcal M_m(z):=\left\{\sum_{j\in J_m} M^2_{(\psi_j,\varphi_j)}(z) \right\}^{1/2}
\end{align*}
and
\begin{align*}
\widetilde T_{mq}:=\widetilde {(T_m)}_q\qquad 
\end{align*}
for ease of notation.

\section{The Range $p\le q$}\label{sec:main1}

In this section we prove Theorem \ref{thm:main1}. The implication (a)$\implies$(b) holds by \eqref{eq:normalized} and \eqref{eq:vanish}, and the implication (c)$\implies$(a) is trivial. So, the only nontrivial implication is  
\begin{align*}
\text{(b)$\implies$ (c)}.
\end{align*}
The rest of the section is devoted to the proof of this implication, to which the key is the gauge control asserted in the next proposition.

\begin{prop}\label{prop:gauge-control}
Let $0< q\le \infty$ and $\{(\psi_j,\varphi_j)\}_{j\in J}$ be a canonical collection of inducing pairs. 
Then the following assertions hold:
\begin{enumerate}
\item[(a)] If $\widetilde T_q\in L^\infty$, 
 \begin{enumerate} 
 \item[(a1)] then $\mathcal M\in L^\infty$;
 \item[(a2)] and if $q<\infty$, then $\mathcal M_0\in C_0$;
 \end{enumerate}
\item[(b)] If, in addition, $\widetilde T_q\in C_0$
\begin{enumerate}
\item[(b1)] and if $q<\infty$, then $\mathcal M\in C_0$; 
\item[(b2)] and if $q=\infty$, then $\mathcal M_1\in C_0$.
\end{enumerate}
\end{enumerate}
\end{prop}

We will provide two independent proofs for Proposition \ref{prop:gauge-control}. Both will be based on the symbol rigidity provided by Proposition \ref{prop:affine}. So, given a collection of pairwise distinct affine symbols, we devise the two independent approaches described below. 
\begin{enumerate}
\item[$\bullet$] (Slope-class Approach) As observed in the proof of Lemma $\ref{lem:rigidity}$, each summand is separated from the full sum by a carefully chosen differential operator. Here, we employ a new differential operator to isolate each slope-class operator and establish slope-class rigidity to complete the proof. This approach is analytic in nature.
\item[$\bullet$] (Gram-matrix Approach)  We have by \eqref{eq:M-equality} and \eqref{eq:duality2}
\begin{align}\label{eq:M-upper}
\mathcal M(z)&\le  \norm{G_{\beta^\w}}|Tk_\w (z)|e^{-\frac{\alpha}{2}|z|^2} \qquad \bigl( \w=\bigl(\varphi_j(z)\bigr)_{j\in J}\bigr)
\end{align}
where $\mathcal M$ is the function specified in \eqref{eq:gaugenorm}. In light of this, we aim to derive a uniform estimate for $\norm{G_{\beta^\w}}$ by exploiting the affine property of the symbols. Overall, this approach is largely algebraic.
\end{enumerate}
\medskip

\bigskip
\noindent
\begin{center}
{\bf $\bullet$ Slope-class Approach}
\end{center}
\smallskip

To begin with, we first establish the rigidity for slope-class operators.

\begin{lem}\label{lem:slope-class}
Let $\{(\psi_j,\varphi_j)\}_{j\in J}$ be a canonical collection of inducing pairs with symbols as in {\rm (N)}. Let 
$$
a_1=\dots=a_N=:a
$$
and suppose that there are polynomials $Q_1, \dots, Q_N$ of degree $d$ and a function $\rho:\C\to [0,\infty)$ such that 
\begin{align}\label{eq:sumcond}
\left|\sum_{j=1}^N \psi_j(z)Q_j(z)e^{\alpha w\varphi_j(z)}\right| \leq \rho(w)(1+|z|+|w|)^d e^{\frac{\alpha}{2}(|z|^2+|w|^2)}
\end{align}
for all $z,w\in\C$. Then the following assertions hold:
\begin{enumerate}
\item[(a)] If $\rho\in L^\infty$, then $\mathcal M\in L^\infty$. 
\item[(b)] If $a\ne0$ and $\lim_{|w|\to\infty}\rho(w)=0$, then $\mathcal M\in C_0$. 
\end{enumerate}
\end{lem}

\begin{proof}
Fix an arbitrary $z\in\C$ and put  
\begin{align}\label{eq:Rzw}
R_z(w):= \sum_{j=1}^N  \psi_j(z)Q_j(z)e^{\alpha b_j (w+\overline{az})}.
\end{align}
We have by assumption \eqref{eq:sumcond}
\begin{align*}
|R_z(w)|
&= \left|e^{-\alpha az(w+\overline{az})}\sum_{j=1}^N \psi_j(z)Q_j(z)e^{\alpha(w +\overline{az})(az+b_j)}\right|\\
&\lesssim \rho(w+\overline{az})(1+|z|+|w+\overline{az}|)^d |e^{-\alpha az(w+\overline{az})}|e^{\frac{\alpha}{2}( |z|^2 +|w+\overline{az}|^2)} \\
&\lesssim \rho(w+\overline{az})(1+|z|+|w|)^d e^{\frac{\alpha}{2}[ (1-|a|^2)|z|^2+|w|^2]}
\end{align*}
for all $z, w\in\C$.
In particular, 
\begin{align}\label{eq:Rzzeta-est}
|R_z(\zeta)|\lesssim \rho(\zeta+\overline{az})(1+|z|)^d e^{\frac{\alpha}{2}(1-|a|^2)|z|^2}
\end{align}
for all $z, \zeta\in \C$ with $|\zeta|\le 1$.

Meanwhile, since $b_1, \dots, b_N$ are pairwise distinct by canonicality, we can find a small $\epsilon\in \bigl(0, \frac{1}{N}\bigr)$ such that $e^{\alpha\epsilon b_1}, \dots, e^{\alpha\epsilon b_N}$ are also pairwise distinct. So, arguing as in the proof of Proposition \ref{prop:weight-rigid}, we obtain by \eqref{eq:Rzw}
\begin{align}\label{eq:Vandermonde}
\begin{pmatrix}  \psi_1(z)Q_1(z)e^{\alpha\overline{az}b_1} \\
\psi_2(z)Q_2(z)e^{\alpha\overline{az}b_2}\\
 \vdots  \\
 \psi_N(z)Q_N(z)e^{\alpha\overline{az}b_N} \\
 \end{pmatrix}
= V^{-1}
 \begin{pmatrix} R_z(0) \\
 R_z(\epsilon)\\
 \vdots  \\
R_z\bigl((N-1)\epsilon\bigr) \\
 \end{pmatrix}
\end{align}
where $V$ is the Vandermonde matrix specified in \eqref{eq:V}.
In addition, setting 
$$
\omega(z):=\max_{0\le m \le N-1} \rho(m\epsilon+\overline{az}),
$$
we obtain by \eqref{eq:Rzzeta-est}
\begin{align}\label{eq:Rzest}
|R_z(m\epsilon)|\lesssim \omega(z) (1+|z|)^d e^{\frac{\alpha}{2}(1-|a|^2)|z|^2},\qquad m=0,\dots, N-1,
\end{align}
with implicit constant independent of $z$.

Now, fix $j$. Since 
$$
\left|e^{\alpha\overline{az}b_j}\right|= e^{\frac{\alpha}{2}(|az+b_j|^2-|az|^2-|b_j|^2)},
$$
we have 
\begin{align*}
M_{(\psi_j, \varphi_j)}(z)
&=|\psi_j(z)|e^{\frac{\alpha}{2} (|az+b_j|^2-|z|^2)}\\
&=\frac{|\psi_j(z)Q_j(z)e^{\alpha \overline{az}b_j}|} {|Q_j(z)|}e^{-\frac{\alpha}{2}(1-|a|^2)|z|^2}e^{\frac{\alpha}{2}|b_j|^2}
\end{align*}
whenever $Q_j$ is non-vanishing at $z$.
In conjunction with this, we see from \eqref{eq:Rzest}
and \eqref{eq:Vandermonde} that
\begin{align*}
|\psi_j(z)Q_j(z)e^{\alpha\overline{az}b_j}|\lesssim\omega(z)(1+|z|)^d e^{\frac{\alpha}{2}(1-|a|^2)|z|^2},
\end{align*}
and hence
\begin{equation*}
M_{(\psi_j,\varphi_j)}(z)\lesssim\omega(z)\frac{(1+|z|)^d}{|Q_j(z)|} \approx \omega(z)\quad\text{as}\ \ |z|\to\infty;
\end{equation*}
for the second estimate, recall that $Q_j$ is a polynomial of degree $d$. 
Since this holds for every $j$, we obtain
\begin{align}\label{eq:M-omega}
\mathcal M(z)\lesssim\omega(z)
\end{align}
for all sufficiently large $|z|$.
If $\rho\in L^\infty$, then $\omega$ is essentially bounded, and hence
$\mathcal M\in L^\infty$ by continuity on bounded sets. Thus (a) holds. 
If, in addition, $a\ne0$, then
$$
|m\epsilon+\overline{az}|\ge |a||z|-m\epsilon \to\infty\quad\text{as}\ \ |z|\to\infty
$$
for each $m=0, \dots, N-1$. Consequently, if $\rho(z)\to 0$ as $|z|\to \infty$, then so does $\omega(z)$.
Thus (b) holds by \eqref{eq:M-omega}. The proof is complete.
\end{proof}

In the proof below, we continue using the notation introduced in the proof of Lemma \ref{lem:rigidity}.

\begin{proof} [\bf First proof of Proposition \ref{prop:gauge-control}]
Noting
$$
F_z(w)=\sum_{j=1}^N \psi_j(z) e^{\alpha w \varphi_j(z)}= e^{\frac{\alpha}{2}\abs{w}^2} T k_{\overline w}(z),
$$
we obtain by \eqref{eq:point-eval}
\begin{align*}
|F_z(w)| \le e^{\frac{\alpha}{2}(\abs{z}^2 +\abs{w}^2)}\widetilde T_q (\overline w)
\end{align*}
for all $z, w\in\C$. Motivated by this, we put
\begin{align*}
\rho(w):=\sup_{|\zeta|\le1}\widetilde T_q (\overline w+\zeta).
\end{align*}

First, we prove (a). So, assume $\widetilde T_q\in L^\infty$. Note $\rho\in L^\infty$. 
By Proposition \ref{prop:affine}, all the symbols are affine. Fix $k$ and set the corresponding slope class 
$$
S_k:=\{j\in J: a_j=a_k\};
$$
we are using the notation from {\rm (N)}. Denote by $|S_k|$ the number of indices in $S_k$. 
We will find a collection $\{Q_j\}_{j\in S_k}$ of polynomials of degree $N-|S_k|$ such that
\begin{align}\label{eq:goal}
\begin{aligned}
&\left|\sum_{j\in S_k} \psi_j(z) Q_j(z) e^{\alpha w\varphi_j(z)}\right|\\
&\qquad\qquad \le C \rho(w) (1+|z|+|w|)^{N - |S_k|} e^{\frac\alpha2(|z|^2+|w|^2)}
\end{aligned}
\end{align}
for some constant $C>0$ independent of $z, w\in\C$. With this granted, we may apply Lemma~\ref{lem:slope-class}(a) to each of the
slope classes to conclude (a1). When $S_k=J$, note from \eqref{eq:Tbdd} that \eqref{eq:goal} holds with $Q_j\equiv1$ and $d=0$.

We now turn to the proof of \eqref{eq:goal} for $S_k\ne J$.
Given $z\in \C$, put
$$
D_k := \prod_{\ell\not\in S_k} \bigl(\partial_w -\alpha \varphi_\ell(z) \bigr)
$$
and 
$$
Q_j(z):=\alpha^{N-|S_k|}\prod_{\ell\not\in S_k} \bigl(\varphi_j(z)-\varphi_\ell(z)\bigr)\qquad (j\in S_k).
$$
Note that $D_k$ is a differential operator of order $N-|S_k|$ and each $Q_j$ is a polynomial of degree $N-|S_k|$. Moreover, 
since $D_k$ annihilates $e^{\alpha w\varphi_j(z)}$ for $j\notin S_k$, we see that
\begin{align*}
D_k F_z(w)=\sum_{j=1}^N \psi_j(z)  D_k e^{\alpha w \varphi_j(z)}=\sum_{j\in S_k} \psi_j(z) Q_j(z) e^{\alpha w\varphi_j(z)}.
\end{align*}
Note 
$$
\sup_{j\in J} |\varphi_j(z)|\lesssim 1+|z|.
$$
Also, for $|\zeta-w|=\frac{1}{1+|w|}\le 1$, we have by \eqref{eq:Tbdd}
\begin{align*}
 |F_z(\zeta)|\le   e^{\frac{\alpha}{2}(\abs{z}^2 +\abs{\zeta}^2)}\widetilde T_q \bigl(\overline \zeta\bigr) \lesssim e^{\frac{\alpha}{2}(\abs{z}^2 +\abs{w}^2)}\rho(w);
\end{align*} 
we used $|\zeta|^2 \le |w|^2 +3$ in the last estimate. Thus, proceeding as in the proof of \eqref{eq:LkFest} with $c_1 = c_2= \frac{\alpha}{2}$, we obtain \eqref{eq:goal}, as desired.
This completes the proof of (a1).

For the proof of (a2), let $q<\infty$ and consider $j\in J_0$ (hence $\varphi_j\equiv b_j$). 
Note from (a1) and Proposition \ref{prop:single-wco}(b) that each summand of $T_1$, and hence $T_1$ itself as well, is bounded from $\F^q$ into itself. Thus $\widetilde T_{1q}\in L^\infty$. Now, since both $\widetilde T_q$ and $\widetilde T_{1q}$ belong to $L^\infty$, so does $\widetilde T_{0q}$.
Accordingly, $\psi_j\in \mathcal F^q$ by Proposition \ref{prop:weight-rigid}
and hence by \eqref{eq:vanish-infty}
$$
|\psi_j(z)|e^{-\frac{\alpha}{2}|z|^2}\to 0\quad\text{as}\ \ |z|\to \infty,
$$
or, said differently, $M_{(\psi_j, b_j)}\in C_0$. This proves (a2).

Finally, if $\widetilde T_q\in C_0$, then Lemma~\ref{lem:slope-class}(b) yields
$$
M_{(\psi_j,\varphi_j)}\in C_0
$$
for each $j\in J_1$. From this and (a2), we see that (b) holds. 
 The proof is complete.
\end{proof}

\bigskip
\noindent
\begin{center}
{\bf $\bullet$ Gram-matrix Approach}
\end{center}
\smallskip

For $\delta>0$, we say that ${\bf w}=(w_j)_{j\in J}\in\C^N$ is a $\delta$-{\em separated} point if 
\begin{align*}
|w_i-w_j|\ge\delta\qquad (i\ne j).
\end{align*}
When the exact value of $\delta$ is not essential, we will simply say ``separated" instead of ``$\delta$-separated". 
This notion of separation naturally extends to sequences, which we will use later in Lemma \ref{lem:separated}. 

Recall that $G_{\beta^\w}$ denotes the Gram matrix associated with the dual basis $\beta^\w$ of $k_\w$; see Section \ref{sec:Gram}. 

\begin{thm}\label{thm:Gram}
Let $\delta>0$. Then $G_{\beta^\w}:\C^N\to\C^N$ is bounded uniformly in $\delta$-separated points ${\bf w}\in \C^N$.
\end{thm}

\begin{proof}
Fix an arbitrary $\delta$-separated point ${\bf w}=(w_j)_{j\in J}\in\C^N$. Put
\begin{align*}
p_j(z):=\prod_{i\ne j}\frac{z-w_i}{w_j-w_i}\quad\text{and}\quad h_j:=p_jk_{w_j}
\end{align*}
for $j\in J$. The reproducing property yields 
\begin{equation*}
\ip{\Pi h_j}{k_{w_i}}=\ip{h_j}{k_{w_i}}=\delta_{ij}\qquad (i, j\in J)
\end{equation*}
where $\Pi:\F^2 \to \operatorname{span}k_\w$ is the orthogonal projection. Thus we obtain 
$$
\beta^\w=(\Pi h_j)_{j\in J}
$$
by the uniqueness of the dual basis.
For each $j$, we note by the $\delta$-separation condition
\begin{align*}
|p_j(w_j+\zeta)|\le \left(1+\frac{|\zeta|}{\delta}\right)^{N-1}
\end{align*}
and thus
\begin{align*}
\frac{\pi}{\alpha} \triplenorm{h_j}_2^2 
&= \int_{\C}|p_j(z)|^2e^{-\alpha|z-w_j|^2}\,dA(z)\quad\text{by}\ \ \eqref{eq:kernelsize}\\
&= \int_\C |p_j(w_j+\zeta)|^2 e^{-\alpha|\zeta|^2}\, dA(\zeta)\\
&\le \int_{\C}
\left(1+\frac{|\zeta|}{\delta}\right)^{2N-2}
e^{-\alpha|\zeta|^2}\,dA(\zeta)<\infty.
\end{align*}
Accordingly, we obtain
\begin{align*}
\norm{\beta^\w_{j}} \le \triplenorm{h_j}_2 \le C\qquad(j\in J)
\end{align*}
for some constant $C=C(\alpha, N,\delta)>0$. Now, we obtain by \eqref{eq:eigen-max}
\begin{align*}
\norm{G_{\beta^\w}}\le \operatorname{tr}(G_{\beta^\w}) =\sum_{j=1}^N \norm{\beta^\w_j}^2 \le N C^2,
\end{align*}
as desired. The proof is complete.
\end{proof}

 In what follows we use the notation 
\begin{align}\label{eq:Gram-dual-uniform}
B_\delta:=\sup \norm{G_{\beta^\w}}<\infty\qquad (\delta>0)
\end{align}
where the supremum is taken over all $\delta$-separated points $\w\in\C^N$. The finiteness of $B_\delta$ comes from Theorem \ref{thm:Gram}.

\begin{proof}[\bf Second Proof of Proposition \ref{prop:gauge-control}]
By Proposition \ref{prop:affine}, symbols are all affine and pairwise distinct. Thus we can find  $\delta>0$ such that
$\w:=\bigl(\varphi_j(z)\bigr)_{j\in J}$ is a $\delta$-separated point for $|z|$ sufficiently large, say $|z|\ge R$.
\\
Fix $z$ with $|z|\ge R$. We have by \eqref{eq:M-upper} and \eqref{eq:Gram-dual-uniform}
\begin{align*}
\mathcal M(z)\le  B_\delta|Tk_\w (z)|e^{-\frac{\alpha}{2}|z|^2}.
\end{align*}
Moreover, by \eqref{eq:point-eval},
\begin{align*}
|Tk_\w (z)|^2e^{-\alpha|z|^2}
=\sum_{j=1}^N |Tk_{\varphi_j(z)}(z)|^2 e^{-\alpha|z|^2}
\le \sum_{j=1}^N \widetilde T_q \bigl(\varphi_j(z)\bigr).
\end{align*}
Combining these observations yields
\begin{align}\label{eq:Gram-M-est}
\mathcal M(z) \le B_\delta \left\{\sum_{j=1}^N \widetilde T^2_q \bigl(\varphi_j(z)\bigr)\right\}^{1/2},
\end{align}
and hence
\begin{align*}
\mathcal M(z)\le B_\delta\sqrt{N}\bigl\|\widetilde T_q\bigr\|_{L^\infty},\qquad |z|\ge R.
\end{align*}
Since gauge functions are already bounded for $|z|\le R$ by continuity, we deduce from the above that (a1) holds.  
The proof for (a2) is the same as in the first proof of Proposition \ref{prop:gauge-control}. 

Next, assume $\widetilde T_q\in C_0$. We now use the notation from {\rm (N)}.
Note that $T_1$ is bounded from $\F^q$ into itself by (a1) and Proposition \ref{prop:single-wco}(b). Thus $\widetilde T_{1q}$, and hence $\widetilde T_{0q}$ as well, belongs to $L^\infty$. 
Accordingly, by Proposition \ref{prop:weight-rigid}, $T_0$ is a finite-rank compact operator from $\F^q$ into itself. 
In particular, $\widetilde T_{0q}\in C_0$ by \eqref{eq:vanish}.
Now, since both $\widetilde T_q$ and $\widetilde T_{0q}$ belong to $C_0$, so does $\widetilde T_{1q}$. 

Meanwhile, applying the argument leading to \eqref{eq:Gram-M-est} to the reduced sum $T_1$, we obtain
\begin{align}\label{eq:reducedgauge}
\mathcal M_1(z) \le B_\delta \left\{\sum_{j\in J_1} \widetilde T^2_{1q} \bigl(\varphi_j(z)\bigr)\right\}^{1/2},\qquad |z|\ge R.
\end{align}
Now, since  $\widetilde T_{1q}\in C_0$ and $|\varphi_j(z)|\to\infty$ as $|z|\to\infty$ for each $j\in J_1$, 
we conclude $\mathcal M_1\in C_0$. 
This completes the proof.
\end{proof}

\bigskip
\noindent
\begin{center}
{\bf $\bullet$ Conclusion}
\end{center}
\smallskip

We are now ready to prove Theorem \ref{thm:main1}. Note that the implications (a)$\implies$(b) and (c)$\implies$(a) are trivial. So, it remains to show that (b) implies (c). 
In the proof below we use the notation from {\rm (N)}.

\begin{proof}[\bf Proof of (b)$\implies$(c) in Theorem \ref{thm:main1}]

Assume (b). By Proposition \ref{prop:affine}, all the symbols are affine. 
Note from Propositions \ref{prop:gauge-control} and \ref{prop:single-wco}(b) that $T_1:\F^p\to\F^q$ is bounded (compact, resp.) and thus $\widetilde T_{1q}\in L^\infty$.
Now, since $\widetilde T_q\in L^\infty$ by assumption, we have $\widetilde T_{0q}\in L^\infty$.
Accordingly, for each $j\in J_0$, we have $\psi_j\in \F^q$ by Proposition \ref{prop:weight-rigid} and hence $W_{\psi_j,b_j}:\F^p\to \F^q$ is compact by Proposition \ref{prop:single-wco}(a). Thus (c) holds, as desired. The proof is complete. 
\end{proof}

\section{The Reverse Range $0<q<p\le\infty$}\label{sec:main2}

In this section we prove Theorem \ref{thm:main2}. The implications (e)$\implies$(b)$\implies$(a) are trivial and the equivalence 
(d)$\iff$ (e) holds by Proposition \ref{prop:single-wco}(a), (c). We will complete the proof by establishing the implications 
\begin{align}\label{eq:implication}
\text{(a)$\implies$ (c) $\implies$ (d)}.
\end{align}

We first prove the second implication of the above, which is relatively short. Throughout the section we put 
$$
s:=\frac{q}{1-\lambda}
$$
for short. Recall $\lambda=\frac{q}{p}\in [0,1)$ with $\lambda=0$ when $p=\infty$.

\begin{proof}[\bf Proof of (c)$\implies$(d) in Theorem \ref{thm:main2}]
Assume that
\begin{align}\label{eq:tildeq-C0}
\widetilde T_q\in L^s.
\end{align}
H\"older's Inequality, applied to the right hand side of \eqref{eq:submean-transform} (with $r=1$), gives
\begin{align*}
\widetilde T^s_q(w)\lesssim\int_{D_1(w)}\widetilde T^s_q(\zeta)\,dA(\zeta).
\end{align*}
The right-hand side is bounded by $\|\widetilde T_q\|_{L^s}^s$ and tends to $0$ as $|w|\to\infty$.
By the continuity property established in \eqref{eq:transform-continuity}, we conclude $\widetilde T_q\in C_0$.
Thus, by Theorem \ref{thm:main1}, each summand is bounded from $\F^q$ into itself. We now use the notation from {\rm (N)}.

By Proposition 2.1(a), each constant-symbol summand is bounded from $\mathcal F^p$ to $\mathcal F^q$.
Moreover, its kernel-norm $q$-transform being a translate of a Gaussian function, we have 
\begin{align}\label{eq:T0-Ls}
\widetilde T_{0q}\in L^s
\end{align}
and hence $\widetilde T_{1q}\in L^s$ by \eqref{eq:tildeq-C0}.

Meanwhile, proceeding as in the second proof of Proposition \ref{prop:gauge-control}, we can find $\delta>0$ and $R>0$ for which  \eqref{eq:reducedgauge} holds. 
Since any two norms are equivalent on a finite dimensional space, we obtain 
\begin{align*}
\mathcal M^s_1(z) \lesssim \sum_{j\in J_1} \widetilde T^s_{1q}(a_jz+b_j), \qquad |z|\ge R,
\end{align*}
so that
\begin{align*}
\int_{|z|\ge R} \mathcal M^s_1(z)\,dA(z)
&\lesssim \sum_{j\in J_1}\int_{\C}\widetilde T^s_{1q}(a_jz+b_j)\,dA(z)\\
&= \left(\sum_{j\in J_1}|a_j|^{-2}\right)\bigl\|\widetilde T_{1q}\bigr\|_{L^s}^s<\infty.
\end{align*}
Since gauge functions are already bounded for $|z|\le R$ by continuity, this yields $\mathcal M_1\in L^s$.
Proposition 2.1(c) now implies that each summand of $T_1$ is bounded from $\F^p$ to $\F^q$, as desired.
Together with the constant-symbol summands already considered, this completes the proof.
\end{proof}

The rest of the section is devoted to the proof of the first implication of \eqref{eq:implication}, which is the hardest part.
Due to the length of the proof, we first outline the overall strategy before presenting the technical details. 

\begin{enumerate}
\item[$\bullet$] (Overall Scheme) The proof is based on a discretized gauge estimate. After removing the constant-symbol summands, we consider the values of the remaining affine symbols on a sufficiently fine lattice outside a fixed disk. The uniform Gram-matrix estimates from Theorem \ref{thm:Gram} provide stable dual systems, while Lemmas 5.4 and 5.5 below yield a synthesis estimate and local control of the combined gauge \(\mathcal M\). Combining these estimates with Khinchine's Inequality, bounded overlap, and sequence-space duality, we obtain
$$
\mathcal M\in L^s.
$$
Together with Young's Convolution Inequality, this will lead us to the desired conclusion.
\end{enumerate}
\medskip

We use the notation $\ell^p$ for the standard sequence space over the set of all positive integers.

\begin{lem}\label{lem:separated}
Let $0<p\le \infty$ and $\epsilon>0$. Let $u=(u_\ell)_{\ell=1}^\infty$ be an $\epsilon$-separated sequence in $\C$. For $c=(c_\ell)_{\ell=1}^\infty\in \ell^p$, put 
$$
f:=\sum_{\ell=1}^\infty c_\ell k_{u_\ell}.
$$
Then $f\in\F^p$ and 
$$
\triplenorm{f}_p\le  C \norm{c}_{\ell^p}
$$
for some constant $C=C(p,\epsilon)>0$.
\end{lem}

\begin{proof}
Let $0< p\le 1$. In this case we have 
$$
|f|^p\le \sum_{\ell=1}^\infty |c_\ell|^p |k_{u_\ell}|^p
$$
Thus, integrating both sides of the above over $\C$ against the relevant Gaussian measure, we obtain by \eqref{eq:kernelsize}
$$
\triplenorm{f}_p^p\le \norm{c}_{\ell^p}^p.
$$
Next, consider the case $p=\infty$. In this case, we note from \eqref{eq:kernelsize}
\begin{align}\label{eq:supnorm}
\triplenorm{f}_\infty\le \norm{c}_{\ell^\infty}\left(\sup_{z\in\C} \sum_{\ell=1}^\infty e^{-\frac{\alpha}{2} |z-u_\ell|^2}\right).
\end{align}
To see that the above supremum is finite, fix $z$ and divide the plane into the annuli
\begin{align*}
A_n(z)
=\left\{\zeta\in\C:\frac{n\epsilon}{2}\le |\zeta-z|<\frac{(n+1)\epsilon}{2}\right\}
\end{align*}
for integers $n\ge 0$. Since the sequence $u$ is $\epsilon$-separated, a standard packing
argument shows that
\begin{align*}
x_n:=\sharp\{\ell: u_\ell\in A_n(z)\}\le C_1( n+1)
\end{align*}
for some constant $C_1=C_1(\epsilon)>0$. Hence
\begin{align*}
\begin{aligned}
\sum_{\ell=1}^\infty e^{-\frac{\alpha}{2}|z-u_\ell|^2}
&\le C_1\sum_{n=0}^\infty (n+1) e^{-\frac{\alpha}{8} \epsilon^2 n^2}<\infty.
\end{aligned}
\end{align*}
This, together with \eqref{eq:supnorm}, proves the lemma for the case $p=\infty$. 

Now, having proved the lemma for the endpoint cases $p=1$ and $p=\infty$, we may now conclude the lemma for the case $1<p<\infty$ by interpolation; see \cite[Section 2.4]{Zhu}. 
The proof is complete. 
\end{proof}

We now introduce further notation. 
For
$\xi=(\xi_j)_{j\in J}$ and $\w=(w_j)_{j\in J}\in \C^N$ with $w_j$'s pairwise distinct, define 
\begin{align}\label{eq:Aw-def}
\Lambda_\w^\xi:=\Bigl(\ip{\beta^\w_i}{k_{\xi_j}}\Bigr)=\Bigl(   \beta^\w_i (\xi_j) e^{-\frac{\alpha}{2}|\xi_j|^2}       \Bigr)\in \Mat.
\end{align}
Let $U_\w^\xi\in\Mat$ be the diagonal unitary matrix whose $j$-th diagonal entry is $\overline\eta_j$, where
\begin{align*}
\eta_j:=\frac{ \left| \ip {k_{\xi_j}   }  { k_{w_j}} \right|} {\ip {k_{\xi_j}   }  { k_{w_j}}}.
\end{align*}
Note $\Lambda_\w^\w=U_\w^\w=I_N$; recall that $I_N$ denotes the identity matrix.

The following lemma shows that $\Lambda_\w^\xi U_\w^\xi$ remains uniformly close to $I_N$ when $\xi$ is close to a separated point $\bf w$.
Recall that $B_\delta$ denotes the constant specified in \eqref{eq:Gram-dual-uniform}.

\begin{lem}\label{lem:uniform-close}
Let $\delta>0$. Then 
\begin{align*}
\norm {\Lambda_\w^\xi U_\w^\xi-I_N}\le \sqrt{\alpha N B_\delta}|\xi-\w|
\end{align*}
for all $\xi\in\C^N$ and $\delta$-separated points ${\bf w}\in\C^N$. 
\end{lem}

\begin{proof}
Consider arbitrary $\xi=(\xi_j)_{j\in J}$ and a $\delta$-separated point $\w=(w_j)_{j\in J}$ in $\C^N$. 
Let $\gamma_j$ be the $j$-th column of $\Lambda_\w^\xi$.
Then the $j$-th column of $\Lambda_\w^\xi U_\w^\xi-I_N$ is 
$$
\overline{\eta_j} \gamma_j-{\bf e}_j
$$
whose $i$-th entry  is  
\begin{align*}
\ip{\beta^\w_i}{\eta_jk_{\xi_j}-k_{w_j}}
&=\ip{E_{\beta^\w} {\bf e}_i}{\eta_jk_{\xi_j}-k_{w_j}}\\
&=\ip{{\bf e}_i}{E^*_{\beta^\w} [\eta_jk_{\xi_j}-k_{w_j}]}_{\C^N},
\end{align*}
where $E_{\beta^\w}: \C^N\to \operatorname{span} k_\w$  denotes the linear operator defined in Section \ref{sec:Gram}. Note from \eqref{eq:positive}
$$
\norm{G_{\beta^\w}}=\|E^*_{\beta^\w} E_{\beta^\w}\| =\|E^*_{\beta^\w}\|^2
$$
where $\|E^*_{\beta^\w}\|$ denotes the operator norm. 
Accordingly, 
\begin{align}\label{eq:column-norm}
|\overline{\eta_j} \gamma_j-{\bf e}_j|= \bigl|E_{\beta^\w}^* [\eta_jk_{\xi_j}-k_{w_j}]\bigr|\le \sqrt{\norm{G_{\beta^\w}}} \triplenorm{\eta_jk_{\xi_j}-k_{w_j}}_2
\end{align}
for each $j$. In conjunction with this, we note  
\begin{align*}
\triplenorm{\eta_j k_{\xi_j}-k_{w_j}}_2^2
&=\triplenorm{ k_{\xi_j}}_2^2+ \triplenorm{k_{w_j}}_2^2 -2\re\ip{\eta_j k_{\xi_j}}{k_{w_j}}\\
&= 2(1-\abs{\ip{k_{\xi_j}}{k_{w_j}}})\\
&=2\bigl(1-e^{-\frac{\alpha}{2}|\xi_j-w_j|^2}\bigr)\\
&\le \alpha |\xi_j-w_j|^2
\end{align*}
for each $j$. 
Combining this with \eqref{eq:column-norm} yields  
\begin{align*}
\begin{aligned}
\bigl|\overline{\eta_j} \gamma_j-{\bf e}_j\bigr|\le \sqrt{\alpha \norm{G_{\beta^\w}}}|\xi-\w|
\end{aligned}
\end{align*}
for each $j$. Now, since the operator norm of a matrix is dominated by its Euclidean norm, we obtain 
\begin{align*}
\|\Lambda_\w^\xi U_\w^\xi-I_N\|
\le \left\{\sum_{j=1}^N \bigl|\overline{\eta_j} \gamma_j-{\bf e}_j\bigr|^2\right\}^{1/2} 
\le\sqrt{\alpha N \norm{G_{\beta^\w}}}|\xi-\w|
\end{align*}
for all $\xi, \w\in \C^N$. Now, further assuming that $\w$ is $\delta$-separated, we conclude the lemma by \eqref{eq:Gram-dual-uniform}.
The proof is complete.
\end{proof}

For the subsequent corollary, we recall the relation
\begin{align}\label{eq:RR}
\inf_{|{\bf z}|=1,\, {\bf z}\in \C^N}\norm {L{\bf z}}= \sigma_{\min}(L),\qquad L\in \Mat,
\end{align}
which is an immediate consequence of the Rayleigh-Ritz Theorem. Here, $\sigma_{\min}$ denotes the smallest singular value.

\begin{cor}\label{cor:sigma_min}
Let $\delta>0$ and put $r_\delta:=(2\sqrt{\alpha N B_\delta})^{-1}$. Then
\begin{equation*}
\inf_{|{\bf z}|=1,\, {\bf z}\in \C^N}\norm {\Lambda_\w^\xi{\bf z}}\ge \frac12 \quad \text{whenever} \quad |\xi-\w|\le r_\delta
\end{equation*}
for all $\xi\in\C^N$ and $\delta$-separated points ${\bf w}\in\C^N$. 
\end{cor}

\begin{proof}
Fix $\xi\in\C^N$ and let ${\bf w}\in\C^N$ be a $\delta$-separated point.
Put $X_\w^\xi :=\Lambda_\w^\xi U_\w^\xi$ for short. 
By Lemma \ref{lem:uniform-close} we have
\begin{align*}
\|X_\w^\xi -I_N\| \le\frac12 \quad \text{whenever} \quad |\xi-\w|\le r_\delta.
\end{align*}
Accordingly, for $|\xi-\w|\le r_\delta$, we have
\begin{align*}
\abs{X_\w^\xi {\bf z}}\ge \abs{\bf z}- \abs{(X_\w^\xi -I_N){\bf z}}\ge \frac 12 |{\bf z}|
\end{align*}
for all ${\bf z}\in\C^N$ and thus by \eqref{eq:RR}
$$
\sigma_{\min}\bigl(X_\w^\xi\bigr)=\inf_{|{\bf z}|=1} \abs{X_\w^\xi  {\bf z}}\ge \frac12.
$$
Since singular values are invariant under right multiplication by a unitary matrix, we conclude the corollary by \eqref{eq:RR} once again.
The proof is complete.
\end{proof}
\bigskip

\subsection*{\quad\bf Technical Preparation (TP)} \label{sec:prep}

We now pause to go through the necessary technical preparation toward
Theorem~\ref{thm:main2}. Unless otherwise emphasized, the assumptions
and notation introduced here will be used tacitly for the rest of the
section.
Since $T: \mathcal F^p\to \mathcal F^q$ is bounded by assumption, we note $\widetilde T_q\in L^\infty$.
Accordingly, by Theorem \ref{thm:main1},  summands are individually bounded from $\F^q$ into itself. 
Thus, by the proof of the implication (c)$\implies$(d), the kernel-norm $q$-transforms of constant-symbol summands all belong to $L^s$.
Thus, after subtracting the constant-symbol summands and relabeling the remaining sum, we may assume by  \eqref{eq:affine}
$$
\varphi_j(z)=a_jz+b_j\quad\text{with}\quad 0<|a_j|\le 1, \qquad j\in J;
$$
we are using the notation from {\rm (N)}. Put  
$$
{\bf a}:=(a_j)_{j\in J}\in\C^N
$$
for the rest of the section. Let 
$$
\Phi(z):=\bigl(\varphi_j(z)\bigr)_{j\in J}\in\C^N.
$$
Since symbols are all affine and pairwise distinct, we can easily find $\delta>0$ and $R>0$ such that 
$\Phi(z)$ is a $\delta$-separated point when $|z|\ge R$.

Now, using the number $r_\delta$ specified in Corollary \ref{cor:sigma_min}, fix $\tau>0$ with 
\begin{align}\label{eq:tau-choice}
\frac{\tau|{\bf a}|}{\sqrt2}<r_\delta
\end{align}
and pick an enumeration $(z_\ell)_{\ell=1}^\infty$ of the standard square $\tau$-lattice 
$$
\left\{\tau(m+n\sqrt{-1}): m, n\ \text{are integers}\right\}.
$$
Let $Q_\ell$ be the half-open standard square of side length $\tau$ centered at $z_\ell$. Let $\Gamma$ be the set of all indices $\ell$ satisfying
\begin{align}\label{eq:retain}
|z_\ell|\ge R+\frac{\tau}{\sqrt2}
\end{align}
The cell $Q_\ell$ is said to be {\em retained} if $\ell\in \Gamma$ and {\em omitted} otherwise. By \eqref{eq:retain}, retained cells all 
contained in $\C\setminus D_R(0)$, while the omitted cells are all contained in the closed disk $\overline D_{R+\sqrt2 \tau}(0)$. Set 
$$
 {\bf w_\ell}:= (w_{j\ell})_{j\in J}\quad \text{where}\quad  w_{j\ell} :=\varphi_j(z_\ell)
$$
so that 
$$
\w_\ell=\Phi(z_\ell).
$$
Note from \eqref{eq:retain} that $\w_\ell$ is $\delta$-separated for each $\ell\in\Gamma$. Hence 
\begin{align}\label{eq:Gell-dual-uniform}
\sup_{\ell\in\Gamma}\norm{G_{\beta^{\w_\ell}}}\le B_\delta<\infty
\end{align}
by \eqref{eq:Gram-dual-uniform}.

\bigskip
With the technical preparation complete, we now resume our progress toward Theorem \ref{thm:main2}.
In what follows we denote by $\ell^p_\Gamma$ the linear subspace of $\ell^p$ consisting of all $\Gamma$-supported sequences in $\ell^p$.

\begin{lem}\label{lem:S}
Let $\Gamma$ be as in {\rm (TP)} and let $0<p\le \infty$. 
Let $u=(u_\ell)_{\ell=1}^\infty\in\ell^p_\Gamma$ be a finitely supported sequence. 
Then 
\begin{align}\label{eq:S}
\triplenorm{\sum_{\ell\in\Gamma} u_\ell \beta^{\w_\ell}_j}_p\le  C \norm{u}_{\ell^p_\Gamma}\qquad (j\in J)
\end{align}
for some constant $C>0$ independent of $j$ and $u$.
\end{lem}

\begin{proof}

For each $i\in J$, since
\begin{align*}
w_{i\ell}=\varphi_i(z_\ell)=a_i z_\ell+b_i
\end{align*}
and the lattice mesh is $\tau$, we have
\begin{align*}
|w_{i\ell }-w_{i m}|=|a_i|\,|z_\ell-z_m|\ge|a_i|\tau \qquad(\ell\ne m).
\end{align*}
So, setting 
$$
\epsilon:=\min_{i\in J} |a_i|\tau >0,
$$
we see that $(w_{i\ell})_{\ell\in\Gamma}$ is an $\epsilon$-separated sequence for each $i$. 

Fix $j\in J$ and $\ell\in \Gamma$. Let $d_{i j}^{(\ell)}$ be the $ij$-entry of $G_{\beta^{\w_\ell}}$. Then we have by \eqref{eq:inverse} 
\begin{align*}
\beta^{\w_\ell}_j=\sum_{i=1}^N d_{i j}^{(\ell)}k_{w_{i\ell}}.
\end{align*}
We also have by \eqref{eq:Gell-dual-uniform}
\begin{align*}
\Bigl|d_{ij}^{(\ell)}\Bigr|\le B_\delta<\infty.
\end{align*}
Now, since $(w_{i\ell})_{\ell\in\Gamma}$ is an $\epsilon$-separated sequence, we obtain by Lemma \ref{lem:separated}
\begin{align*}
\triplenorm{\sum_{\ell\in\Gamma}u_\ell d_{ij}^{(\ell)}k_{w_{i\ell}}}_p
\lesssim\left\|\bigl(u_\ell d_{ij}^{(\ell)}\bigr)_\ell\right\|_{\ell^p_\Gamma} 
\le B_\delta \|u\|_{\ell^p_\Gamma}
\end{align*}
with implicit constant independent of $i$ and $j$.

Finally, since $u$ is finitely supported, we have for each $j$
\begin{align*}
\sum_{\ell\in\Gamma}u_\ell \beta^{\w_\ell}_j=\sum_{i=1}^N\left(\sum_{\ell\in\Gamma}u_\ell d_{ij}^{(\ell)}k_{w_{i\ell}}\right)
\end{align*}
and hence conclude \eqref{eq:S} by the two types of triangle inequalities depending on the cases $0<p<1$ and $1\le p\le\infty$. The proof is complete. 
\end{proof}

Note from \eqref{eq:M-equality}
\begin{align*}
\mathcal M(z_\ell)=|T\beta^{\w_\ell}(z_\ell)|e^{-\frac{\alpha}{2}|z_\ell|^2}
\end{align*}
for all $\ell\in\Gamma$; recall $T\beta^{\w_\ell}= (T\beta^{\w_\ell}_j)_{j\in J}$.   When the lattice mesh is sufficiently small, the exact duality stated above\,--\,which holds  at the lattice center of a retained cell\,--\,persists throughout the cell as shown in the next lemma. 

\begin{lem}\label{lem:LP}
Let $\{(\psi_j,\varphi_j)\}_{j\in J}$ be a canonical collection of inducing pairs as in {\rm (TP)}.
Let $Q_\ell$ be the cells as in  {\rm (TP)}. Then 
$$
\mathcal M(z)\le 2 |T\beta^{\w_\ell}(z)|e^{-\frac{\alpha}{2}|z|^2},\qquad z\in Q_\ell,
$$
for all $\ell\in\Gamma$.
\end{lem}

\begin{proof}
Fix $\ell\in\Gamma$ and let $z\in Q_\ell$. 
Put $\w:=\Phi(z)$. Recalling $Q_\ell\subset \C\setminus D_R(0)$ by \eqref{eq:retain}, we note that $\w$ is $\delta$-separated and hence $\beta^\w$ is well-defined.
Recall from \eqref{eq:M-equality}
$$
T\beta^{\w} (z)  =\bigl(\psi_j(z) e^{\frac{\alpha}{2} |\varphi_j(z)|^2}\bigr)_{j\in J}.
$$
Also, recall from \eqref{eq:Aw-def}
\begin{align*}
\Lambda_{\w_\ell}^{\w}=\Bigl(   \beta^{\w_\ell}_i \bigl(\varphi_j(z)\bigr) e^{-\frac{\alpha}{2}|\varphi_j(z)|^2} \Bigr)\in\Mat.
\end{align*}
So, the $i$-th entry of $\Lambda_{\w_\ell}^{\w}[T\beta_{\w} (z)]$  is 
\begin{align*}
\sum_{j=1}^N  \beta^{\w_\ell}_i \bigl(\varphi_j(z)\bigr) e^{-\frac{\alpha}{2}|\varphi_j(z)|^2} \cdot \psi_j(z) e^{\frac{\alpha}{2} |\varphi_j(z)|^2}
= T\beta^{\w_\ell}_i(z).
\end{align*}
In other words, 
\begin{align}\label{eq:Y-A-c}
\Lambda_{\w_\ell}^{\w} [T\beta^\w(z)]=  T\beta^{\w_\ell}(z).
\end{align}
Meanwhile, since $\varphi_j(z)-\varphi_j(z_\ell)= a_j (z-z_\ell)$ and $z\in Q_\ell$, we have 
$$
|\w-\w_\ell|= |{\bf a}||z-z_\ell|\le \frac{\tau|{\bf a}|}{\sqrt2}\le r_\delta
$$
where the last inequality holds by the choice \eqref{eq:tau-choice} of $\tau$. 
Now, since each $\w_\ell$ is a $\delta$-separated point, we see from \eqref{eq:Y-A-c}, Corollary \ref{cor:sigma_min}, and \eqref{eq:M-equality} that
\begin{align*}
2 |T\beta^{\w_\ell}(z)|e^{-\frac{\alpha}{2}|z|^2}\ge |T\beta^\w(z)|e^{-\frac{\alpha}{2}|z|^2} =\mathcal M(z),
\end{align*}
as required. The proof is complete. 
\end{proof}

For the proof of Proposition \ref{prop:reverse-integrability} below, we recall the well-known Khinchine's Inequality.
For a positive integer $n$, the $n$-th Rademacher function $R_n$ is defined by
$$
R_n(t):=\text{sgn}\{\sin(2^n \pi t ) \}, \qquad 0 \le t \le 1.
$$
For a proof of the following Khinchine's Inequality, see \cite[Appendix C]{G}.

\begin{lem}[Khinchine's Inequality]
Given $0 < q <\infty$, there is a constant $C=C(q)>0$ such that
$$
C^{-1}\left( \sum_{n=1}^\infty |c_n|^2 \right)^{q/2}
\le \int_0^1 \left|\sum_{n=1}^\infty c_n R_n(t) \right|^q \,dt
\le C \left(\sum_{n=1}^\infty |c_n|^2 \right)^{q/2}
$$
	for all $(c_n)_{n=1}^\infty \in \ell^2$.
\end{lem}

We are now ready to prove the main part of Theorem \ref{thm:main2}.

\begin{prop}
\label{prop:reverse-integrability}
Let $0<q<p\le\infty$ and $\{(\psi_j,\varphi_j)\}_{j\in J}$ be a canonical collection of inducing pairs as in {\rm (TP)}.
Suppose that $T:\F^p\to \F^q$ is bounded. Then $\mathcal M\in L^s$ where $\lambda:=\frac{q}{p}$.

\end{prop}

\begin{proof}
We divide the argument into three steps.

\medskip
\noindent
\textbf{Step 1: Randomization.}
Fix $j\in J$. Let $u=(u_\ell)_{\ell=1}^\infty\in\ell^p_\Gamma$ be a finitely supported sequence. For $t\in[0,1]$, put
\begin{align*}
f_t = \sum_{\ell=1}^\infty R_\ell(t)u_\ell \beta^{\w_\ell}_j= \sum_{\ell\in\Gamma} R_\ell(t)u_\ell \beta^{\w_\ell}_j,
\end{align*}
where $R_\ell$ denotes the $\ell$-th Rademacher function. By Lemma~\ref{lem:S} we have 
\begin{align*}
\triplenorm{f_t}_p
\lesssim
\|R_\ell(t) u\|_{\ell^p_\Gamma}= \|u\|_{\ell^p_\Gamma}
\end{align*}
and hence 
\begin{align*}
\triplenorm{Tf_t}_q \lesssim \norm{T}_{\F^p\to\F^q} \|u\|_{\ell^p_\Gamma}
\end{align*}
uniformly in $t$. Integrating with respect to $t\in[0,1]$ and using Fubini's theorem gives
\begin{align*}
\int_{\C}\left\{\int_0^1\left|\sum_{\ell\in\Gamma}R_\ell(t)u_\ell T\beta^{\w_\ell}_j(z)\right|^q dt\right\}\, e^{-\frac{q\alpha}{2}|z|^2}\,dA(z)
\lesssim
\|u\|_{\ell^p_\Gamma}^q
\end{align*}
and, moreover, the inner integral is comparable to 
\begin{align*}
\left\{\sum_{\ell\in\Gamma}|u_\ell|^2|T\beta^{\w_\ell}_j(z)|^2\right\}^{q/2}
\end{align*}
by Khinchine's Inequality. Consequently,
\begin{equation}\label{eq:Khinchine-square}
\int_{\C}\left\{\sum_{\ell\in\Gamma}|u_\ell|^2|T\beta^{\w_\ell}_j(z)|^2e^{-\alpha|z|^2}\right\}^{q/2}\,dA(z)
\lesssim \|u\|_{\ell^p_\Gamma}^q;
\end{equation}
with implicit constant independent of $j$ and $u$.

\medskip
\noindent
\textbf{Step 2: Localization.}
Fix $\epsilon>0$ and define
\begin{align*}
Q_\ell^*
:=\{\zeta\in\C:\operatorname{dist}(\zeta,Q_\ell)<\epsilon\}
=\bigcup_{z\in Q_\ell} D_\epsilon(z),\qquad \ell\in\Gamma.
\end{align*}
Since $(z_\ell)_{\ell\in \Gamma}$ is a part of a square lattice, the family $(Q_\ell^*)_{\ell\in\Gamma}$ has bounded overlap, i.e.,
\begin{equation}\label{eq:bounded-overlap}
\sum_{\ell\in\Gamma} \rchi_\ell\le n_0\quad\text{on}\quad \C
\end{equation}
for some positive integer $n_0=n_0(\tau, \epsilon)>0$. Here, $\rchi_\ell$ is the characteristic function of $Q_\ell^*$. 

Since any two norms on an $N$-dimensional space are equivalent, we have by Lemma~\ref{lem:LP} 
\begin{align*}
\mathcal M^q(z) 
\lesssim \left\{ \sum_{j=1}^N|T\beta^{\w_\ell}_j(z)|^2e^{-\alpha|z|^2}\right\}^{q/2}
\approx \sum_{j=1}^N |T\beta^{\w_\ell}_j(z)|^q e^{-\frac{q\alpha}{2}|z|^2}
\end{align*}
for all $\ell\in\Gamma$ and $z\in Q_\ell$. So, applying \eqref{eq:local} on $D_\epsilon(z)\subset Q_\ell^*$ to each entire function
$T\beta^{\w_\ell}_j$, we obtain
\begin{align*}
\mathcal M^q(z)\lesssim\sum_{j=1}^N\int_{Q_\ell^*}|T\beta^{\w_\ell}_j(w)|^qe^{-\frac{q\alpha}{2}|w|^2}\,dA(w)
\end{align*}
with hidden constant independent of $\ell\in\Gamma$ and $z\in Q_\ell$. Thus, we obtain 
\begin{align}\label{eq:S-localized}
S_{Q_\ell}:
&=\sup_{z\in Q_\ell}\mathcal M^q(z)\nonumber\\
&\lesssim\sum_{j=1}^N\int_{Q_\ell^*}|T\beta^{\w_\ell}_j(w)|^qe^{-\frac{q\alpha}{2}|w|^2}\,dA(w)
\end{align}
for all $\ell\in\Gamma$. Multiply \eqref{eq:S-localized} by $|u_\ell|^q$ and sum over all $\ell\in\Gamma$. We obtain 
\begin{align*}
\begin{aligned}
\sum_{\ell\in\Gamma}|u_\ell|^q S_{Q_\ell}
&\lesssim\sum_{j=1}^N\int_{\C}\sum_{\ell\in\Gamma} \mathbf \rchi_\ell(w)|u_\ell|^q|T\beta^{\w_\ell}_j(w)|^qe^{-\frac{q\alpha}{2}|w|^2}\,dA(w).
\end{aligned}
\end{align*}
For fixed $w$ and $j$, put
\begin{align*}
c_\ell := |u_\ell| |T\beta^{\w_\ell}_j(w)| e^{-\frac{\alpha}{2}|w|^2}.
\end{align*}
One may verify
\begin{align*}
\left\{\sum_{\ell\in\Gamma}\rchi_\ell (w) c_\ell^q\right\}^{2/q }\le C\sum_{\ell\in\Gamma} c_\ell^2
\end{align*}
where 
\begin{align*}
C=C(n_0,q):=
\begin{cases}
n_0^{\frac{2}{q}-1},&0<q\le2,\\
1,&q>2.
\end{cases}
\end{align*}
In fact, one may apply H\"older's Inequality and \eqref{eq:bounded-overlap} for $q\le2$, and the triangle inequality for  $q> 2$.
Therefore
\begin{align*}
\begin{aligned}
\sum_{\ell\in\Gamma} |u_\ell|^q S_{Q_\ell}
&\lesssim \sum_{j=1}^N\int_{\C}\left\{\sum_{\ell\in\Gamma}|u_\ell|^2|T\beta^{\w_\ell}_j(w)|^2 e^{-\alpha|w|^2} \right\}^{q/2}\,dA(w).
\end{aligned}
\end{align*}
Applying \eqref{eq:Khinchine-square} to each term of the right hand side of the above, we conclude 
\begin{equation}\label{eq:sequence-master}
\sum_{\ell\in\Gamma}|u_\ell|^q  S_{Q_\ell}\lesssim\|u\|_{\ell^p_\Gamma}^q
\end{equation}
with implicit constant independent of $u$.  

\medskip
\noindent
\textbf{Step 3: Gauge Integrability via Duality.}
Assume first that $p<\infty$. Put $v:=(v_\ell)_{\ell=1}^\infty$ where $v_\ell:=|u_\ell|^q$.
Note
\begin{align*}
\|u\|_{\ell^p_\Gamma}^q = \left( \sum_{\ell\in\Gamma} v_\ell^{1/\lambda} \right)^{\lambda} = \|v\|_{\ell^{1/\lambda}_\Gamma}.
\end{align*}
Thus \eqref{eq:sequence-master} becomes
\begin{align*}
\sum_{\ell\in\Gamma} v_\ell S_{Q_\ell}\lesssim\|v\|_{\ell^{1/\lambda}_\Gamma}
\end{align*}
with implicit constant independent of $v$. Moreover, since $u\in\ell^p_\Gamma$ is finitely supported and arbitrary, the above remains valid for any finitely supported $v\in \ell^{1/\lambda}_\Gamma$ with nonnegative terms. Since $p<\infty$, finitely supported sequences form a dense subset of $\ell^{1/\lambda}_\Gamma$.  
So, by the duality $\bigl(\ell^{1/\lambda}_\Gamma\bigr)^*=\ell^{\frac{1}{1-\lambda}}_\Gamma$ (recall $1<\frac{1}{\lambda}=\frac{p}{q}<\infty$), we obtain 
\begin{equation}\label{eq:S-dual}
(S_{Q_\ell})_{\ell\in\Gamma} \in \ell^{\frac{1}{1-\lambda}}_\Gamma.
\end{equation}
Now, since the retained cells are pairwise disjoint and have area $\tau^2$, we obtain 
\begin{align*}
\begin{aligned}
\int_{\bigcup_{\ell\in\Gamma} Q_\ell}\mathcal M^s\,dA
=\sum_{\ell\in\Gamma}\int_{Q_\ell}\mathcal M^s\,dA
\le\tau^2\sum_{\ell\in\Gamma} S^{\frac{1}{1-\lambda}}_{Q_\ell}<\infty.
\end{aligned}
\end{align*}
On the other hand, recalling that omitted cells are all contained in $\overline D_{R+\sqrt2 \tau}(0)$, we have by Proposition~\ref{prop:gauge-control}(a), 
\begin{align*}
\int_{\bigcup_{\ell\notin\Gamma} Q_\ell} \mathcal M^s\,dA\le \pi(R+\sqrt2 \tau)^2\norm{\mathcal M}^s_{L^\infty} <\infty. 
\end{align*}
Consequently, we conclude
\begin{align}\label{eq:M-integrable}
\mathcal M\in L^s,
\end{align}
as desired. 

It remains to consider $p=\infty$. In
\eqref{eq:sequence-master}, take $u$ to be the characteristic
sequence of an arbitrary finite set $F\subset \Gamma$.
Then
\begin{align*}
\sum_{\ell\in F}S_{Q_\ell} \lesssim 1
\end{align*}
with hidden constant independent of $F$. Thus \eqref{eq:S-dual} remains valid with $\lambda=0$ and so does \eqref{eq:M-integrable}. 
This completes the proof.
\end{proof}

We are now ready to complete the proof of Theorem \ref{thm:main2} by establishing the first implication in \eqref{eq:implication}.

\begin{proof}[\bf Proof of (a)$\implies$(c) in Theorem \ref{thm:main2}]

Assume that $T: \mathcal F^p\to \mathcal F^q$ is bounded. Note $\widetilde T_q\in L^\infty$ by \eqref{eq:normalized}. As noted in (TP), we may assume that all the symbols are nonconstant. 
In the proof below we simplify the notation by setting
\begin{align*}
W_j:=W_{\psi_j,\varphi_j},\quad M_j:=M_{(\psi_j,\varphi_j)}, \quad\widetilde W_{jq} :=\widetilde{(W_j)}_q 
\end{align*}
for $j\in J$. Note $\mathcal M\in L^s$ by Proposition \ref{prop:reverse-integrability} and hence $M_j\in L^s$ for each $j\in J$.

Fix $j\in J$. Note from \eqref{eq:kernelsize}
\begin{align*}
\widetilde W^q_{jq}(w) =\frac{q\alpha}{2\pi} \int_{\C}M^q_j(z) e^{-\frac{q\alpha}{2} |a_jz +b_j-w|^2}\,dA(z).
\end{align*}
Thus, setting 
\begin{align*}
G_q(\zeta):=\frac{q\alpha}{2\pi} e^{-\frac{q\alpha}{2}|\zeta|^2}\quad\text{and}\quad H_j(\zeta):=|a_j|^{-2}M^q_j\left(\frac{\zeta-b_j}{a_j}\right),
\end{align*}
we see by a linear change of variables 
$$
\widetilde W^q_{jq}=G_q*H_j.
$$
As $\norm{G_q}_{L^1}=1$ and $\frac{1}{1-\lambda}\ge1$, Young's Convolution Inequality yields
\begin{align*}
\norm{\widetilde W_{jq}}_{L^s}^q
&=\norm{G_q*H_j}_{L^{\frac{1}{1-\lambda}}}\\
&\le\norm{H_j}_{L^{\frac{1}{1-\lambda}}}\\
&=|a_j|^{-2\lambda}\norm{M_j}_{L^s}^q
\end{align*}
so that
$$
\norm{\widetilde W_{jq}}_{L^s} \le |a_j|^{-\frac{2}{p}} \norm{M_j}_{L^s} <\infty,
$$
where $|a_j|^{-\frac{2}{p}}=1$ if $p=\infty$. Since this holds for each $j\in J$, we conclude $\widetilde T_q\in L^s$, as desired.
This completes the proof. 
\end{proof}

Now that Theorems \ref{thm:main1} and \ref{thm:main2} have been established, a retrospective look back at Proposition \ref{prop:single-wco} for a single operator reveals that one can employ kernel-norm transforms in place of gauge functions. We conclude the paper with a remark that provides direct details in this direction.

\begin{rem}\label{remark}
Given an inducing pair $(\psi, \varphi)$ with $\varphi(z)=az+b$ ($|a|\le 1$), put $W:=W_{\psi,\varphi}$ and $M:=M_{(\psi, \varphi)}$ for simplicity. 
Since
\begin{align}\label{eq:equality}
|Wk_w(z)| e^{-\frac{\alpha}{2}|z|^2} = M(z) e^{-\frac{\alpha}{2}|az+b-w|^2},
\end{align}
we have 
\begin{align}\label{eq:Wqint}
\widetilde W^q_q(w)=\frac{q\alpha}{2\pi}\int_{\C}M^q(z)e^{-\frac{q\alpha}{2}|az+b-w|^2}\,dA(z)
\end{align}
for $0<q<\infty$ and 
\begin{align}\label{eq:Winfty}
\widetilde W_\infty(w)&=\sup_{z\in\C}M(z)e^{-\frac{\alpha}{2}|az+b-w|^2}.
\end{align}
Note from \eqref{eq:point-eval}
\begin{align}\label{eq:MW}
M(z)\le \widetilde W_q(az+b)
\end{align}
for $0<q\le\infty$.

First, consider the case $a\ne 0$. In this case we have the equivalences 
\begin{align}\label{eq:MWLt}
M\in L^t \iff \widetilde W_q\in L^t\qquad (0<q\le t\le\infty)
\end{align}
and 
\begin{align}\label{eq:MWC0}
M\in C_0
\iff \widetilde W_q\in C_0\qquad (0<q\le\infty).
\end{align}
As a consequence of these equivalences, one may employ kernel-norm transforms in place of gauge functions for the single-operator characterization when $a\ne 0$ in Proposition \ref{prop:single-wco}. 
{\em(Sketch of Proof)}: The equivalence \eqref{eq:MWLt} for $q<t<\infty$ is already mentioned at the end of Section \ref{sec:single-operator}. At the endpoints when $q<\infty$, we have
\begin{align*}
\bigl\|\widetilde W_q\bigr\|_q = \|M\|_q\quad\text{by \eqref{eq:Wqint}}
\end{align*}
and 
\begin{align*}
|a|^{\frac{2}{q}}\bigl\|\widetilde W_q\bigr\|_{L^\infty}\le \|M\|_{L^\infty}\le \|\widetilde W_q\|_{L^\infty}\quad\text{by \eqref{eq:Wqint} and \eqref{eq:MW}},
\end{align*}
which remains valid even when $q=\infty$ by \eqref{eq:Winfty}.
Meanwhile, since the integral in \eqref{eq:Wqint} is a Gaussian convolution after the change of variables $\zeta=az+b$ and convolution with an $L^1$-function preserves $C_0$, we obtain 
\eqref{eq:MWC0} for $q<\infty$. The same holds for $q=\infty$ by \eqref{eq:equality} and \eqref{eq:MW} due to the Gaussian decay. $\Box$.

Now, suppose $a=0$ so that $\varphi(z)\equiv b$. In this case, we have 
\begin{align*}
\widetilde W_q (w) = \triplenorm{\psi}_q e^{\frac{\alpha}{2}|b|^2 } e^{-\frac{\alpha}{2}|w-b|^2},
\end{align*}
and so 
\begin{align*}
\bigl\|\widetilde W_q\bigr\|_{L^t} = c_t e^{\frac{\alpha}{2}|b|^2 }\triplenorm{\psi}_q\qquad (0<t\le \infty)
\end{align*}
where $c_t:=\left(\frac{2\pi}{t\alpha}\right)^{1/t}$ for $t<\infty$ and $c_\infty=1$. Consequently, the membership condition $\psi\in\F^q$ is equivalent to each of the conditions 
$$
\widetilde W_q\in L^\infty, \quad \widetilde W_q\in C_0, \quad \widetilde W_q\in L^t
$$
for any $0<t<\infty$. This is consistent with the rank-one characterization in Proposition \ref{prop:single-wco}(a). 
\end{rem}

\medskip
\noindent
{\quad\bf Acknowledgement.} As noted in the Introduction, the symbol rigidity established in Section~\ref{sec:pre-rigidity} serves as a common starting point for the proofs of our main results.
We are grateful to Song-Ying Li for inspiring the idea behind this symbol rigidity.

\medskip
\textbf{AI Statement.}
The authors used artificial-intelligence tools for language editing,
\LaTeX{} formatting, and limited assistance in checking local mathematical
arguments. The proof strategy and mathematical development are entirely the
authors' own. All mathematical statements, proofs, and conclusions were
independently verified by the authors, who take full responsibility for
the content of the paper.

\end{document}